\documentclass[10pt,onecolumn,draftcls]{IEEEtran}

\pdfoutput = 1

\usepackage[cmex10]{amsmath}
\usepackage{dsfont}
\usepackage{algorithm}
\usepackage{algorithmic}
\usepackage[pdftex]{graphicx}

\DeclareMathOperator*{\argmin}{arg\,min}

\usepackage{lineno}

\begin{document}
%\linenumbersep 5pt\relax
%\linenumbers

\title{
Convex optimization problem prototyping
for image reconstruction in computed tomography
with the Chambolle-Pock algorithm 
}

\author{Emil Y. Sidky,
\thanks{Emil Y. Sidky and Xiaochuan Pan are with the Department of Radiology,
University of Chicago, 5841 S. Maryland Ave., Chicago, IL 60637, USA.
e-mail: \{sidky,xpan\}@uchicago.edu.}%
Jakob~H.~J{\o}rgensen,
\thanks{Jakob H. J{\o}rgensen is with the Department of Informa\-tics and Mathematical Modeling,
Technical University of Denmark, Richard Petersens Plads, Building 321, 2800 Kgs. Lyngby, Denmark.
e-mail: jakj@imm.dtu.dk.}%
and Xiaochuan Pan}% <-this % stops a space

\maketitle

\begin{abstract}
The primal-dual optimization
algorithm developed in Chambolle and Pock (CP), 2011 is applied to various
convex optimization problems of interest in computed tomography (CT) image reconstruction.
This algorithm allows for rapid prototyping of optimization problems for the purpose
of designing iterative image reconstruction algorithms for CT.  The primal-dual
algorithm is briefly summarized in the article, and its potential for prototyping
is demonstrated by explicitly deriving CP algorithm instances for many optimization
problems relevant to CT.  An example application modeling breast CT
with low-intensity X-ray illumination is presented.
\end{abstract}

\section{Introduction}

Optimization-based image reconstruction algorithms for CT have
been investigated heavily recently due to their potential to
allow for reduced scanning effort while maintaining or improving
image quality \cite{McCollough:09,PanIP:09}.
Such methods have been considered for many years,
but within the past five years computational barriers have been
lowered enough such that iterative image reconstruction can be considered
for practical application in CT \cite{ziegler2008iterative}.
The transition to practice
has been taking place alongside further theoretical development
particularly with algorithms based on the sparsity-motivated
$\ell_1$-norm
\cite{li2002accurate,SidkyTV:06,Sidky2008image,Chen2008prior,SidkyPC:10,ritschl2011improved,Defrise:11,Fessler:2011,jakob:2011}.
Despite the recent interest in sparsity,
optimization-based image reconstruction algorithm
development continues to proceed along many fronts
and there is as of yet no consensus on
a particular optimization problem for the CT system.
In fact, it is beginning to look like the optimization problem,
upon which the iterative image reconstruction algorithms are
based, will themselves be subject to design depending on the
particular properties of each scanner type and imaging task.

Considering the possibility of tailoring optimization problems
to a class of CT scanners, makes design of iterative image
reconstruction algorithms a daunting task. Optimization formulations
generally construct an objective function comprised of a data
fidelity term and possible penalty terms discouraging unphysical
behavior in the reconstructed image, and they possibly include
hard constraints on the image.  The image estimate is arrived
at by extremizing the objective subject to any constraints placed
on the estimate.  The optimization problems for image reconstruction
can take many forms depending on image representation, projection model, and
objective and constraint design. On top of this, it is difficult
to solve many of the optimization problems of interest. A change
in optimization problem formulation can mean many weeks or months
of algorithm development to account for the modification.

Due to this complexity, it would be
quite desirable to have an algorithmic tool to facilitate design
of optimization problems for CT image reconstruction. This tool
would consist of a well-defined set of mechanical steps that
generate a convergent algorithm from a specific optimization problem
for CT image reconstruction. The goal of this tool would be to
allow for rapid prototyping of various optimization formulations;
one could design the optimization problem free of any restrictions imposed
by a lack of an algorithm to solve it.
The resulting algorithm might not be the most efficient solver for the particular
optimization problem, but it would be guaranteed to give the answer.

In this article we consider convex optimization problems for CT image
reconstruction,
including non-smooth objectives, unconstrained and constrained formulations.
One general algorithmic
tool is to use steepest descent or projected steepest descent \cite{Nocedal:06}.
Such algorithms, however, do not address non-smooth objective functions
and they have difficulty with constrained optimization, being applicable for only
simple constraints such as non-negativity.
%The trouble with such algorithms is that accurate solution to many problems
%of interest for CT image reconstruction
%can be prohibitively expensive in terms of computation time.
Another general strategy
involves some form
of evolving quadratic approximation to the objective. The literature on
this flavor of algorithm design is enormous, including non-linear
conjugate gradient methods \cite{Nocedal:06}, parabolic surrogates \cite{Defrise:11,erdogan1999ordered},
and iteratively reweighted
least-squares \cite{green1984iteratively}. For the CT system these strategies
often require quite a bit of know-how due to the very large scale 
and ill-posedness of the imaging model. Once the optimization formulation is established,
however, these quadratic methods provide a good option to gain in efficiency.

One of the main barriers to prototyping alternative optimization problems for CT image
reconstruction is the size of the imaging model; volumes can contain millions of voxels
and the sinogram data can correspondingly consist of millions of X-ray transmission measurements.
For large-scale systems there has been some resurgence of first-order methods
\cite{yin2008bregman,combettes2008proximal,Beck:09,becker2010templates,chambolle2011first,jakob:11}
and recently there has been applications of first-order methods specifically
for optimization-based image reconstruction in CT \cite{jakob:11,Choi:10,joergensen2011toward}.
These methods are interesting because they can be adapted to
a wide range of optimization problems involving non-smooth functions such as
those involving $\ell_1$-based norms.
In particular,
the algorithm that we pursue further in this paper is a first-order
primal-dual algorithm for convex problems by Chambolle and Pock
\cite{chambolle2011first}.  This algorithm goes a long way toward the
goal of optimization problem prototyping, because it covers a very general class of optimization
problems that contain many optimization formulations of interest to the CT community.

For a selection of optimization problems of relevance to CT image reconstruction,
we work through the details of setting up the
Chambolle-Pock algorithm.
We refer to these dedicated algorithms as \emph{algorithm instances}.
Our numerical results demonstrate that the
algorithm instances
achieve the solution of difficult convex optimization problems under challenging
conditions in reasonable time and without parameter tuning.
In Sec. \ref{sec:chambolle} the CP methodology and
algorithm is summarized; in Sec. \ref{sec:CPCT}
various optimization problems for CT image reconstruction are presented along with
their corresponding CP algorithm instance; and Sec. \ref{sec:simulations} shows a limited study
on a breast CT 
simulation that demonstrates the application of the derived CP algorithm instances.

%The main goal is to have
%a framework for optimization formulation prototyping. Thus the emphasis is on
%ease of algorithm-instance derivation, implementation, and paucity of parameters.
%Efficiency is only needed in so much
%as that the solution of the optimization problem can be obtained in "reasonable time".
%Once the desired optimization problem is identified, then attention can be directed
%to improving efficiency.
%Algorithm parameters are a bane for the present purpose. It often occurs
%that an algorithm can achieve convergence in reasonable time, if algorithm parameters
%are chosen appropriately. But it is better, here, to have a program run overnight knowing that
%the solution will be achieved in the morning than to have a program that can achieve
%the solution in ten minutes only if the right algorithm parameter values are
%selected and not at
%all if the wrong parameter values are chosen.
%The Chambolle-Pock (CP) Algorithm 1 appears to suit our purpose. The algorithm does require
%a little more mathematical sophistication than gradient-based methods in that
%familiarity with the Legendre transform, or equivalently convex conjugation,
%is needed.

\section{Summary of the generic Chambolle-Pock algorithm}
\label{sec:chambolle}

The Chambolle-Pock (CP) algorithm \cite{chambolle2011first}
is primal-dual meaning that it solves an optimization
problem simultaneously with its dual.  On the face of it, it would seem to
involve extra work by solving two problems instead of one, but the algorithm
comes with convergence guarantees and solving both problems provides
a robust, non-heuristic convergence check -- the duality gap.

The CP algorithm applies to a general form of the primal minimization:
\begin{linenomath}
\begin{equation}
\label{primal}
\min_{x} \; \left\{ F(Kx) + G(x) \right\} , 
\end{equation}
\end{linenomath}
and a dual maximization:
\begin{linenomath}
\begin{equation}
\label{dual}
\max_{y} \; \left\{ -F^*(y) - G^*(-K^T y) \right\},
\end{equation}
\end{linenomath}
where $x$ and $y$ are finite dimensional vectors in the spaces $X$ and $Y$, respectively;
$K$ is a linear transform from $X$ to $Y$;
$G$ and $F$ are convex, possibly non-smooth, functions mapping the respective
$X$ and $Y$ spaces to non-negative real numbers; and
the superscript ``$^*$'' in the dual maximization problem refers to convex conjugation,
defined in Eqs. (\ref{legendre1}) and (\ref{legendre2}).
We note that the matrix $K$ need not be square; $X$ and $Y$ will in general have
different dimension.
Given a convex function $H$ of a vector $z \in Z$,
its conjugate can be computed by the Legendre transform \cite{rockafellar1970convex}, and
the original function can be recovered by applying conjugation again:
\begin{linenomath}
\begin{align}
\label{legendre1}
H^*(z) & = \max_{z^\prime}  \; \left\{ \langle z, z^\prime  \rangle_Z - H(z^\prime) \right\}, \\
\label{legendre2}
H(z^\prime) & = \max_{z}  \; \left\{ \langle z^\prime,z  \rangle_Z - H^*(z) \right\} .
\end{align}
\end{linenomath}
The notation $\langle \cdot , \cdot \rangle_Z$ refers to the inner product in the vector space $Z$.

Formally, the primal and dual problems are connected
in a generic saddle point optimization problem:
\begin{linenomath}
\begin{equation}
\label{minmax}
\min_{x} \max_{y} \; \left \{ \langle Kx,y \rangle_Y + G(x) - F^*(y) \right\} .
\end{equation}
\end{linenomath}
By performing the maximization over $y$ in Eq. (\ref{minmax}),
using Eq. (\ref{legendre2}) with
$Kx$ associated with $y^\prime$, the primal minimization Eq. (\ref{primal}) is derived.
Similarly, performing the minimization over $x$ in Eq. (\ref{minmax}),
using Eq. (\ref{legendre1}) and
the identity $\langle Kx,y \rangle = \langle x,K^T y \rangle$,  yields the dual
maximization Eq. (\ref{dual}),
where the $T$ superscript denotes matrix transposition.

The minimization problem in Eq. (\ref{primal}), though compact, covers
many minimization problems of interest to tomographic image reconstruction.
Solving the dual problem, Eq. (\ref{dual}), simultaneously allows for assessment
of algorithm convergence. For intermediate estimates $x$ and $y$ of the primal
minimization and the dual maximization, respectively, the primal objective will be
greater than or equal to the dual objective. The difference between these objectives is
referred to as the duality gap, and convergence is achieved when this gap is zero.
Plenty of examples of useful optimization problems for tomographic image reconstruction
will be described in detail in Sec. \ref{sec:CPCT}, but
first we summarize Algorithm 1 from Ref. \cite{chambolle2011first}.

\subsection{Chambolle-Pock: Algorithm 1}

%\begin{algorithm}
%\caption{The power method: yields dominant eigenvector $\hat{f}$ and its singular value.
%$\hat{f}_0$ is an initial image estimate with a norm of 1.}
%\label{alg1}
%\begin{algorithmic}[1]
%\STATE $\hat{f} \gets \hat{f}_0$
%\REPEAT
%\STATE $\vec{f} \gets X^T X \hat{f}$
%\STATE $\hat{f} \gets \vec{f}/|\vec{f}|$
%\UNTIL{no change in $\hat{f}$}
%\STATE $s_\text{max} \gets |X \hat{f}|$
%\end{algorithmic}
%\end{algorithm}

\begin{algorithm}
\caption{Pseudocode for $N$-steps of the basic Chambolle-Pock algorithm.
The constant $L$ is the $\ell_2$-norm of the matrix $K$; $\tau$ and $\sigma$ are non-negative CP
algorithm parameters, which are both set to $1/L$ in the present application; $\theta \in [0,1]$
is another CP algorithm parameter, which is set to 1; and $n$ is the iteration index.
The proximal operators $prox_\sigma$ and $prox_\tau$ are defined in Eq. (\ref{proximal}).}
\label{alg1}
\begin{algorithmic}[1]
\STATE $L \gets \|K\|_2; \; \tau \gets 1/L; \; \sigma \gets 1/L; \; \theta \gets 1; \; n \gets 0$
\STATE initialize $x_0$ and $y_0$ to zero values
\STATE $\bar{x}_0 \gets x_0$
\REPEAT
\STATE $y_{n+1} \gets prox_\sigma[F^*](y_n + \sigma K \bar{x}_n)$ \label{yupdate}
\STATE $x_{n+1} \gets prox_\tau[G](x_n - \tau K^T y_{n+1})$ \label{xupdate}
\STATE $\bar{x}_{n+1} \gets x_{n+1} + \theta(x_{n+1} - x_n)$
\STATE $n \gets n+1$
\UNTIL{$n \ge N$}
\end{algorithmic}
\end{algorithm}

The CP algorithm simultaneously solves Eqs. (\ref{primal}) and (\ref{dual}).
As presented in Ref. \cite{chambolle2011first} the algorithm is simple, yet extremely
effective.  We repeat the steps here in Listing \ref{alg1} for completeness,
providing the parameters that we use for all results shown below.
The parameter descriptions are provided in Ref. \cite{chambolle2011first}, but
note that in our usage specified above there are {\it no} free parameters.
This is an extremely important feature for our purpose of optimization prototyping.  
One caveat is that technically the proof of convergence for the CP algorithm
assumes $L^2 \sigma \tau <1$, but in practice we have never encountered a case
where the choice $\sigma=\tau = 1/L$ failed to tend to convergence.
We stress that in Eq. (\ref{xupdate}) the matrix $K^T$ needs to be the transpose
of the matrix $K$; this point can sometimes be confusing because $K$ for imaging
applications is often intended to be an approximation to some continuous operator
such as projection or differentiation and often $K^T$ is taken to mean the approximation
to the continuous operator's adjoint, which may or may not be the matrix transpose
of $K$.
The constant $L$ is the magnitude of the matrix $K$, its largest singular value.
Appendix \ref{sec:power} gives the details on computing $L$ via the power method.
Key to deriving the particular
algorithm instances are the proximal mappings $prox_\sigma[F^*]$ and
$prox_\tau[G]$ (called resolvent operators in Ref. \cite{chambolle2011first}).

The proximal mapping is used to generate a descent direction for
the convex function $H$ and it is obtained by the following minimization:
\begin{linenomath}
\begin{equation}
\label{proximal}
prox_\sigma[H](z) = \argmin_{z^\prime}  \left\{ 
H(z^\prime) + \frac{\| z - z^\prime \|_2^2}{2 \sigma} \right\}
\end{equation}
\end{linenomath}
This operation does admit non-smooth convex functions, but
$H$ does need to be simple enough that the above minimization
can be solved in closed form. For CT applications the ability to handle non-smooth
$F$ and $G$ allows the study of many optimization problems of recent interest, and the
simplicity limitation is not that restrictive as will be seen.
%In case that it is simpler to compute the
%proximal mapping from a function's convex conjugate, Ref. \cite{chambolle:2011}
%points out that Moreau's identity can be employed:
%\begin{equation}
%\label{Moreau}
%prox_\sigma[F^*](y) = y - \sigma prox_{1/\sigma} [F] (y/\sigma).
%\end{equation}

\subsection{The CP algorithm for prototyping of convex optimization problems}
\label{sec:prototyping}

To prototype a particular convex optimization problem for CT image reconstruction
with the CP algorithm, there are five basic steps:\\
(1) Map the optimization problem to the generic minimization problem in Eq. (\ref{primal}).\\
(2) Derive the dual maximization problem, Eq. (\ref{dual}), by computing the convex conjugates of $F$ and $G$
using the Legendre transform Eq. (\ref{legendre1}). \\
(3) Derive the proximal mappings of $F$ and $G$ using Eq. (\ref{proximal}).\\
(4) Substitute the results of (3) into the generic CP algorithm in Listing \ref{alg1} to obtain
a CP algorithm instance.\\
(5) Run the algorithm, monitoring the primal-dual gap for convergence.

As will be seen below, a great variety of constrained and unconstrained optimization
problems can be written in the form of Eq. (\ref{primal}).
Specifically, using the algebra of convex functions \cite{rockafellar1970convex},
that the sum of two convex functions is convex and that the composition of 
a convex function with a linear transform is a convex function,
many interesting optimization formulations
can be put in the form of Eq. (\ref{primal}). 
We will also make use of 
convex functions which are
not smooth -- notably $\ell_1$ based norms and indicator functions $\delta_S(x)$:
\begin{linenomath}
\begin{equation}
\label{indicator}
\delta_S(x) \equiv
\begin{cases}
0 & x \in S \\
\infty & x \notin S
\end{cases},
\end{equation}
\end{linenomath}
where $S$ is a convex set. The indicator function is particularly handy for
imposing constraints.
In computing the convex conjugate and proximal mapping of convex functions, we make much use
of the standard calculus rule for extremization, $\nabla f = 0$, but such computations
are augmented also with geometric reasoning, which may be unfamiliar.
Accordingly, we have included appendices to show some of these computation steps.
With this quick introduction, we are now in a position to derive various
algorithm instances for CT image reconstruction from different convex optimization problems.

\section{Chambolle-Pock algorithm instances for CT}
\label{sec:CPCT}

For this article, we only consider optimization problems involving the linear
imaging model for X-ray projection, where the data are considered as line integrals
over the object's X-ray attenuation coefficient.
Generically, maintaining consistent notation with Ref. \cite{chambolle2011first}, the 
discrete-to-discrete CT system model \cite{Barrett:FIS}
can be written as:
\begin{linenomath}
\begin{equation}
\label{linsys}
Au=g,
\end{equation}
\end{linenomath}
where $A$ is the projection matrix taking an object represented by expansion coefficients
$u$ and generating a set of line-integration values $g$.  This model covers a multitude of
expansion functions and CT configurations, including both 2D fan-beam and 3D cone-beam projection
data models.

A few notes on notation are in order. In the following, we largely avoid indexing of the
various vector spaces in order that the equations and pseudocode listings are brief and clear.
Any of the standard algebraic operations between vectors is to be interpreted in a component-wise manner
unless explicitly stated. Also, an algebraic operation between a scalar and a vector is to be distributed
among all components of the vector; e.g., $1 + v$ adds one to all components of $v$.
For the optimization problems below, we employ three vector spaces: $I$ the space of discrete images
in either 2 or 3 dimensions;
$D$ the space of the CT sinograms (or projection data); and $V$ the
space of spatial-vector-valued image arrays, $V=I^d$ where $d=2$ or 3
for 2D and 3D-space, respectively. For the CT system model Eq. (\ref{linsys}),
$u\in I$, and $g \in D$, but we note that the space $D$ can also include sinograms which
are not consistent with the linear system matrix $A$. The vector space $V$ will be used below for
forming the total variation (TV) semi-norm;
an example of such a vector $v \in V$ is the spatial-gradient of an image $u$.
Although the pixel representation is used, much of the following can be
applied to other image expansion functions.
As we will be making much use of certain indicator functions, we define
two important sets, $Box(a)$ and $Ball(a)$, through their indicator function:
\begin{linenomath}
\begin{equation}
\label{Boxindicator}
\delta_{Box(a)}(x) \equiv
\begin{cases}
0 & \|x\|_\infty \le a \\
\infty & \|x\|_\infty > a
\end{cases},
\end{equation}
\end{linenomath}
and
\begin{linenomath}
\begin{equation}
\label{Ballindicator}
\delta_{Ball(a)}(x) \equiv
\begin{cases}
0 & \|x\|_2 \le a \\
\infty & \|x\|_2 > a
\end{cases}.
\end{equation}
\end{linenomath}
Recall that the $\| \cdot \|_\infty$ norm selects the largest component of the argument,
thus $Box(a)$ comprises vectors with no component larger than $a$ (in 2D $Box(a)$ is a square
centered on the origin with width $2a$). We also employ $\mathbf{0}_X$
and $\mathbf{1}_X$ to mean
a vector from the space $X$ with all components set to 0 and 1, respectively.

\subsection{Image reconstruction by least-squares}

Perhaps the simplest optimization method for performing image reconstruction is to minimize
the the quadratic data error function. We present this familiar case in order to gain some
experience with the mechanics of deriving CP algorithm instances, and because the quadratic data error
term will play a role in other optimization problems below. The primal problem of interest is:
\begin{linenomath}
\begin{equation}
\label{LSprimal}
\min_u \; \frac{1}{2} \|Au-g\|_2^2.
\end{equation}
\end{linenomath}
To derive the CP algorithm instance,  we make
the following mechanical associations with the primal problem Eq. (\ref{primal}):
\begin{linenomath}
\begin{align}
\label{LSassoc}
F(y) & = \frac{1}{2} \| y - g \|_2^2, \\
G(x) & = 0,\\
x  & = u, \; \; y = Au \\
K &= A .
\end{align}
\end{linenomath}
Applying Eq. \ref{legendre1}, we obtain the convex conjugates of $F$ and $G$:
\begin{linenomath}
\begin{align}
\label{LSconjugates}
F^*(p) & = \frac{1}{2} \|p\|^2_2 + \langle p, g \rangle_D , \notag \\
G^*(q) & = \delta_{\mathbf{0}_I}(q),
\end{align}
\end{linenomath}
where $p\in D$ and $q \in I$.
%, and the $\mathbf{0}_I$ subscript of the indicator function is the set containing only the
%0-data array (all components set to the value 0).
While obtaining $F^*$ in this case involves elementary calculus for extremization of Eq. (\ref{legendre1}),
finding $G^*$ needs some comment for those unfamiliar with convex analysis.
Using the definition of the Legendre transform for $G(x)=0$, we have:
\begin{linenomath}
\begin{equation}
\label{glegendre}
G^*(q) = \max_{x}  \; \langle q, x  \rangle_I.
\end{equation}
\end{linenomath}
There are two possibilities: (1) $q = \mathbf{0}_I$, in which case the maximum value of $\langle q, x  \rangle_I$
is 0, and (2) $q \neq \mathbf{0}_I$, in which case this inner product can increase without bound, resulting
in a maximum value of $\infty$. Putting these two cases together yields the indicator function
in Eq. (\ref{LSconjugates}).
With $F$, $G$, and their conjugates, the optimization problem dual to Eq. (\ref{LSprimal})
can be written down from Eq. (\ref{dual}):
\begin{linenomath}
\begin{equation}
\label{LSdual}
\max_p \; \left\{ - \frac{1}{2} \|p\|_2^2 - \langle p, g \rangle_D - \delta_{\mathbf{0}_I}(- A^T p) \right\}.
\end{equation}
\end{linenomath}
For deriving the CP algorithm instance, it is not strictly necessary to have this dual problem, but
it is useful for evaluating convergence.

The CP algorithm solves
Eqs. (\ref{LSprimal}) and (\ref{LSdual}) simultaneously. In principle, the values of
the primal and dual objective functions provide a test of convergence.
During the iteration
the objective of the primal problem will by greater than the objective of the dual problem,
and when the solutions of the respective problems are reached,
these objectives will be equal. Comparing the duality gap, i.e. the difference
between the primal objective and the dual objective, with 0 thus
provides a test of convergence.
The presence of the indicator function in the dual problem, however, complicates this
test.
Due to the negative sign in front of the indicator, when the argument is not
the zero vector, this term and therefore the whole dual objective is assigned to a value of
$-\infty$. The dual objective achieves a finite, testable value only when the indicator function
attains the value of 0, when $A^T p = \mathbf{0}_I.$ Effectively,
the indicator function becomes a way to write down a constraint in the form of a
convex function, in this case
an equality constraint. The dual optimization problem can thus alternately be
written as a conventional constrained
maximization:
\begin{linenomath}
\begin{equation}
\label{LSdualalt}
\max_p \; \left\{ - \frac{1}{2} \|p\|_2^2 - \langle p, g \rangle_D \right\}
\text{  such that  } A^Tp=\mathbf{0}_I.
\end{equation}
\end{linenomath}
The convergence check is a bit problematic, because the equality constraint will
not likely be strictly satisfied in numerical computation. Instead, we introduce
a conditional primal-dual gap (the difference between the primal and dual objectives
ignoring the indicator function)
given the estimates $u^\prime$ and $p^\prime$:
\begin{linenomath}
\begin{equation}
\label{LSpdgap}
cPD(u^\prime, p^\prime) = \|Au^\prime -g\|^2_2 + 
\frac{1}{2} \| p^\prime\|_2^2 + \langle p^\prime, g \rangle_D,
\end{equation}
\end{linenomath}
and separately monitor $A^T p^\prime$ to see if it is tending to $\mathbf{0}_I$. Note that
the conditional primal-dual gap need not be positive, but it should tend to zero.

To finally attain the CP algorithm instance for image reconstruction by least-squares,
we derive lines \ref{yupdate} and \ref{xupdate} in Alg. \ref{alg1}.
The proximal mapping $prox_\sigma[F^*](y)$, $y \in D$, for this problem results from a
quadratic minimization:
\begin{linenomath}
\begin{align}
\label{LSproxy}
prox_\sigma[F^*](y) &= 
\argmin_{y^\prime} \;  \left\{ \frac{1}{2} \|y^\prime\|_2^2 +
\langle y^\prime, g \rangle_D + \frac{\| y - y^\prime \|_2^2}{2 \sigma} \right\} \\
&= \frac{y- \sigma g}{1+ \sigma}, \notag
\end{align}
\end{linenomath}
and as $G(x)=0$, $x \in I$, the corresponding proximal mapping is
\begin{linenomath}
\begin{equation}
\label{LSproxx}
prox_\tau[G](x) =x.
\end{equation}
\end{linenomath}
Substituting in the arguments from the generic algorithm,
leads to the update steps in Listing  \ref{algls}.
\begin{algorithm}
\caption{Pseudocode for $N$-steps of the least-squares Chambolle-Pock algorithm instance.}
\label{algls}
\begin{algorithmic}[1]
\STATE $L \gets \|A\|_2; \; \tau \gets 1/L; \; \sigma \gets 1/L; \; \theta \gets 1; \; n \gets 0$
\STATE initialize $u_0$ and $p_0$ to zero values
\STATE $\bar{u}_0 \gets u_0$
\REPEAT
\STATE $p_{n+1} \gets (p_n + \sigma( A \bar{u}_n-g))/(1+\sigma)$ \label{plsupdate}
\STATE $u_{n+1} \gets u_n - \tau A^T p_{n+1}$ \label{ulsupdate}
\STATE $\bar{u}_{n+1} \gets u_{n+1} + \theta(u_{n+1} - u_n)$
\STATE $n \gets n+1$
\UNTIL{$n \ge N$}
\end{algorithmic}
\end{algorithm}
The constant $L=\|A\|_2$ is the largest singular value of $A$ (see Appendix \ref{sec:power} 
for details on the power method). Crucial to the implementation of the CP algorithm instance
is that $A^T$ be the exact transpose of $A$, which is a non-trivial matter for tomographic
applications, because the projection matrix $A$ is usually
computed on-the-fly \cite{Siddon:1985,de2004distance,Mueller:07}.
Convergence of the CP algorithm is only guaranteed when $A^T$ is the exact transpose of $A$,
although it may be possible to extend the CP algorithm to mismatched projector/back-projector
pairs by employing the analysis in Ref. \cite{zeng2000unmatched}.

This derivation of the CP least-squares algorithm instance illustrates
the method on a familiar optimization problem, and it provides
a point of comparison with standard algorithms; this
quadratic minimization problem can be solved straight-forwardly
with the basic, linear conjugate gradients (CG) algorithm.
Another important point for this particular algorithm instance, where limited projection
data can lead to an underdetermined system, is that the CP algorithm will yield
a minimizer of the objective $\| A u - g\|^2_2$ which depends on the initial image $u_0$.
In this case, it is recommended to take advantage of the prototyping capability of the
CP framework to augment the optimization problem so that it selects a unique image
independent of initialization.  For example, one often seeks an image closest to either
$\mathbf{0}_I$ or a prior image, which can be formulated by adding a quadratic term $\| u \|^2_2$ or
$\| u - u_{prior} \|^2_2$ with a small combination coefficient.

\subsubsection{Adding in non-negativity constraints}

One of the flexibilities of the CP method becomes apparent in adding bound constraints.
While CG is also flexible tool for dealing with large and small quadratic optimization,
modification to include constraints, such as non-negativity, considerably complicates
the CG algorithm.  For CP, adding in bound constraints is simply a matter of
introducing the appropriate indicator function into the primal problem:
\begin{linenomath}
\begin{equation}
\label{LSPOSprimal}
\min_u \; \left\{ \frac{1}{2} \|Au-g\|_2^2 + \delta_P(u) \right\},
\end{equation}
\end{linenomath}
where the set $P$ is all $u$ with non-negative components.
Again, we make the
mechanical associations with the primal problem Eq. (\ref{primal}):
\begin{linenomath}
\begin{align}
\label{LSPOSassoc}
F(y) & = \frac{1}{2} \| y - g \|_2^2, \\
G(x) & = \delta_P(x),\\
x &= u, \; \; y =Kx ,\\
K  & = A.
\end{align}
\end{linenomath}
The difference from the unconstrained problem is the function $G(x)$.
It turns out that the convex conjugate of $\delta_P(x)$ is:
\begin{linenomath}
\begin{equation}
\label{legpos}
\delta^*_P(x) = \delta_P(-x),
\end{equation}
\end{linenomath}
see Appendix \ref{sec:legendre} for insight on convex conjugate of indicator functions.
Straight substitution of $G^*$ and $F^*$ into Eq. (\ref{dual}), yields
the dual problem:
\begin{linenomath}
\begin{equation}
\label{LSPOSdual}
\max_p \; \left\{ - \frac{1}{2} \|p\|_2^2 - \langle p, g \rangle_D - \delta_P(A^T p) \right\}.
\end{equation}
\end{linenomath}
As a result the conditional primal-dual gap is the same as before. The
difference now is that the constraint checks are that $A^Tp$ and $u$ should
be non-negative.

To derive the algorithm instance, we need the proximal mapping $prox_\tau[G]$, which by
definition is:
\begin{linenomath}
\begin{equation}
\label{deltaprox1}
prox_\tau[\delta_P](x) = \argmin_{x^\prime} \;
\left\{ \delta_P(x^\prime) + \frac{\| x - x^\prime \|^2_2}{2 \tau} \right\}.
\end{equation}
\end{linenomath}
The indicator in the objective prevents consideration of negative components of $x^\prime$.
The $\ell_2$ term can be regarded as a sum over the square difference between components
of $x$ and $x^\prime$; thus the objective is separable and can be
minimized by constructing $x^\prime$ such that
$x^\prime_i = x_i$ when $x_i>0$ and $x^\prime_i = 0$ when $x_i \le 0$. Thus this proximal
mapping becomes a non-negativity thresholding on each component of $x$:
\begin{linenomath}
\begin{equation}
\label{deltaprox2}
[prox_\tau[\delta_P](x)]_i = [pos(x)]_i \equiv
\begin{cases}
x_i & x_i >0 \\
0   & x_i \le 0
\end{cases}.
\end{equation}
\end{linenomath}
Substituting into the generic pseudocode yields Listing \ref{alglspos}.
\begin{algorithm}
\caption{Pseudocode for $N$-steps of the least-squares with non-negativity constraint, CP algorithm
instance.}
\label{alglspos}
\begin{algorithmic}[1]
\STATE $L \gets \|A\|_2; \; \tau \gets 1/L; \; \sigma \gets 1/L; \; \theta \gets 1; \; n \gets 0$
\STATE initialize $u_0$ and $p_0$ to zero values
\STATE $\bar{u}_0 \gets u_0$
\REPEAT
\STATE $p_{n+1} \gets (p_n + \sigma( A \bar{u}_n-g))/(1+\sigma)$ \label{plsposupdate}
\STATE $u_{n+1} \gets pos( u_n - \tau A^T p_{n+1})$ \label{ulsposupdate}
\STATE $\bar{u}_{n+1} \gets u_{n+1} + \theta(u_{n+1} - u_n)$
\STATE $n \gets n+1$
\UNTIL{$n \ge N$}
\end{algorithmic}
\end{algorithm}
Again, we have $L=\|A\|_2$. The indicator function $\delta_P$ leads to
the intuitive modification that non-negativity thresholding is introduced
in line \ref{ulsposupdate} of Listing \ref{alglspos}. In this case the non-negativity constraint in $u$
will be automatically satisfied by all iterates $u_n$.
Upper bound constraints are equally simple to include.

\subsection{Optimization problems based on the Total Variation (TV) semi-norm}

Optimization problems with the TV semi-norm have received much attention
for CT image reconstruction
lately because of their potential to provide high quality images from
sparse view sampling
\cite{SidkyPC:10,Bian:10,Choi:10,ritschl2011improved,han2011algorithm,xia:043706,sidky2011special}.
The TV semi-norm has
been known to be useful for performing edge-preserving regularization,
and recent developments in compressive sensing have sparked even greater
interest in the use of this semi-norm. Algorithm-wise the TV semi-norm
is difficult to handle. Although it is convex, it is not linear,
quadratic or even everywhere-differentiable,
and the lack of differentiability precludes the use of standard
gradient-based optimization algorithms.
%The recent interest in this norm, however, has led to remarkable
%advances in optimization algorithms to the point where we are
%beginning to see recipe-like methods for setting up efficient, large-scale
%solvers for optimization problems involving the TV-norm.
In this sub-section
we go through, in detail, the derivation of a CP algorithm instance for
a TV-regularized least squares data error norm.
We then consider the Kullback-Leibler (KL) data divergence, which is implicitly
employed by many iterative algorithms based on maximum likelihood
expectation maximization (MLEM).  We also consider a data error norm based
on $\ell_1$ which can have some advantage in reducing the impact
of image discretization error, which generally leads to a highly non-uniform
error in the data domain.
Finally, we derive a CP algorithm instance for constrained TV-minimization,
which is mathematically
equivalent to the least-squares-plus-TV
problem \cite{elad2010sparse}, but whose data-error constraint parameter has more
physical meaning than the parameter used in the corresponding unconstrained minimization.
While the previous CP instances solve optimization problems, which can be solved
efficiently by well-known algorithms, the following CP instances are new for the
application of CT image reconstruction.

The optimization problem of interest is
\begin{linenomath}
\begin{equation}
\label{LSTVprimal}
\min_u \; \left\{ \frac{1}{2} \|Au-g\|_2^2 + \lambda \left\| (|\nabla u|) \right\|_1 \right\},
\end{equation}
\end{linenomath}
where the last term, the $\ell_1$-norm of the gradient-magnitude image,
is the isotropic TV semi-norm. The spatial-vector image $\nabla u$ represents a discrete
approximation to the image gradient which is in the vector space $V$,
i.e., the space of spatial-vector-valued image arrays.
The expression $|\nabla u|$ is the gradient-magnitude image, an image array whose
pixel values are the gradient magnitude at the pixel location. Thus, $\nabla u \in V$
and $|\nabla u| \in I$.
Because $\nabla$ is defined in terms of finite differencing,
it is a linear transform from an image array
to a vector-valued image array, the precise form of which
is covered in Appendix \ref{sec:graddiv}. This problem was not explicitly covered
in Ref. \cite{chambolle2011first}, and we fill in the details here.
For this case, matching the primal problem to Eq. (\ref{primal})
is not as obvious as the previous examples.  We recognize in
Eq. (\ref{LSTVprimal}) that both terms involve a linear transform, thus
the whole objective function can be written in the form $F(Kx)$ with the
following assignments:
\begin{linenomath}
\begin{align}
\label{LSTVassoc}
F(y,z) & = F_1(y) + F_2(z), \; \; \;
F_1(y) = \frac{1}{2} \| y - g \|_2^2, \; \; \;
F_2(z) = \lambda \left\|(| z|) \right\|_1, \\
G(x) & = 0,\\
x &= u, \; \; y = Au, \;\; z = \nabla u, \\
K  & = \binom{A}{\nabla} ,
\end{align}
\end{linenomath}
where $u \in I$, $y \in D$, and $z \in V$.
Note that $F(y,z)$ is convex because it is the sum of two convex functions.
Also the linear transform $K$ takes an image vector $x$ and gives a data vector $y$ and
an image gradient vector $z$.
The transpose of $K$, $K^T=(A^T,-div)$,
will produce an image vector from a data vector $y$ and an
image gradient vector $z$:
\begin{linenomath}
\begin{equation}
\label{Ktranspose}
x \leftarrow A^T y - div \,z,
\end{equation}
\end{linenomath}
where we use the same convention as in Ref. \cite{chambolle2011first} that
$-div \equiv \nabla ^T$, see Appendix \ref{sec:graddiv}.

In order to get the convex conjugate of $F$ we need $F_2^*$.
For readers unfamiliar with the Legendre transform of indicator functions
Appendix \ref{sec:legendre} illustrates the transform of some common cases.
By definition,
\begin{linenomath}
\begin{equation}
\label{TVconj}
F_2^*(q) = \max_z \; \left\{ \langle q, z \rangle_V - \lambda \left\| (|z|) \right\|_1 \right\},
\end{equation}
\end{linenomath}
where $q \in V$, like $z$, is a vector-valued image array.
There are two cases to consider: (1) the magnitude image $|q|$ at all pixels
is less than or equal to $\lambda$, i.e. $|q| \in Box(\lambda)$  and
(2) the magnitude image $|q|$ has at least
one pixel greater than $\lambda$, i.e. $|q| \not\in Box(\lambda)$.
It turns out that for the former case the maximization in Eq. (\ref{TVconj})
yields 0, while the latter cause yields $\infty$.
%The function $- \lambda \left\| (|z|) \right\|_1$ in the objective can be visualized
%as the sum of upside down cones at each image pixel. The effect of adding
%the linear term $\langle q, z \rangle_V$ is to tip and distort these cones according
%to the magnitude and direction of $q$ at each pixel. In the former case, all of the cones remain
%pointing ``downwards'' and the objective in Eq. (\ref{TVconj}) is maximized
%by $z=\mathbf{0}_V$, and accordingly takes the value of $0$. In the latter case, the cones
%at pixels where $|q|>\lambda$ are tipped enough that part of them point ``upward'';
%in this direction the objective can increase without bound as a function of $z$.
Putting these two cases together, we have
\begin{linenomath}
\begin{equation}
\label{TVconj2}
F_2^*(q) = \delta_{Box(\lambda)}(|q|).
\end{equation}
\end{linenomath}
The conjugates of $F$ and $G$ are:
\begin{linenomath}
\begin{align}
\label{LSTVconjugates}
F^*(p,q) & = \frac{1}{2} \|p\|^2_2 + \langle p, g \rangle_D + \delta_{Box(\lambda)}(|q|) , \\
G^*(r) & = \delta_{\mathbf{0}_I}(r),
\end{align}
\end{linenomath}
where $p \in D$, $q \in V$, and $r \in I$.

The problem dual to Eq. (\ref{LSTVprimal}) becomes:
\begin{linenomath}
\begin{equation}
\label{LSTVdual}
\max_{p,q} \; \left\{ - \frac{1}{2} \|p\|_2^2 - \langle p, g \rangle_D -
\delta_{Box(\lambda)}(|q|) - \delta_{\mathbf{0}_I}(div\, q- A^T p) \right\}.
\end{equation}
\end{linenomath}
The resulting conditional primal-dual gap is
\begin{linenomath}
\begin{equation}
\label{LSTVpdgap}
cPD(u^\prime, p^\prime,q^\prime) = \frac{1}{2} \|Au^\prime -g\|^2_2 + 
\lambda \left\|( |\nabla u|) \right\|_1 +
\frac{1}{2} \| p^\prime\|_2^2 + \langle p^\prime, g \rangle_P
\end{equation}
\end{linenomath}
with additional constraints $|q^\prime| \in Box(\lambda)$ and $A^T p^\prime - div \, q^\prime=\mathbf{0}_I$.
The final piece needed for putting together the CP algorithm instance for Eq. (\ref{LSTVprimal})
is the proximal mapping:
\begin{linenomath}
\begin{equation}
\label{LSTVprox}
prox_\sigma[F^*](y,z) = \left( \frac{y- \sigma g}{1+ \sigma} ,
\frac{\lambda z}{\max(\lambda \mathbf{1}_I,|z|)} \right).
\end{equation}
\end{linenomath}
The proximal mapping of the data term was covered previously, and that of the TV term is
explained in Appendix \ref{sec:proximals}.
With the necessary pieces in place, the CP algorithm instance for the $\ell_2^2$-TV objective can
be written down in Listing \ref{alglstv}.
\begin{algorithm}
\caption{Pseudocode for $N$-steps of the $\ell_2^2$-TV CP algorithm instance.}
\label{alglstv}
\begin{algorithmic}[1]
\STATE $L \gets \|(A,\nabla)\|_2; \;\tau \gets 1/L; \; \sigma \gets 1/L; \; \theta \gets 1; \; n \gets 0$
\STATE initialize $u_0$, $p_0$, and $q_0$ to zero values
\STATE $\bar{u}_0 \gets u_0$
\REPEAT
\STATE $p_{n+1} \gets (p_n + \sigma( A \bar{u}_n-g))/(1+\sigma)$ \label{plstvupdate}
\STATE $q_{n+1} \gets \lambda (q_n + \sigma \nabla \bar{u}_n )/
\max(\lambda \mathbf{1}_I,|q_n + \sigma \nabla \bar{u}_n |)$ \label{qlstvupdate}
\STATE $u_{n+1} \gets  u_n - \tau A^T p_{n+1} + \tau div \, q_{n+1}$ \label{ulstvupdate}
\STATE $\bar{u}_{n+1} \gets u_{n+1} + \theta(u_{n+1} - u_n)$
\STATE $n \gets n+1$
\UNTIL{$n \ge N$}
\end{algorithmic}
\end{algorithm}
Line  \ref{qlstvupdate}, and the corresponding expression in Eq. (\ref{LSTVprox}),
require some explanation, because the division operation is non-standard
as the numerator is in $V$ and the denominator is in $I$. The effect of this line is to 
threshold the magnitude of the spatial-vectors at each pixel in $q_n + \sigma \nabla \bar{u}_n$
to the value $\lambda$: spatial-vectors larger than $\lambda$ have their magnitude rescaled
to $\lambda$. The resulting thresholded, spatial-vector image is then assigned to $q_{n+1}$.
Recall that $\mathbf{1}_I$ at line \ref{qlstvupdate} is an image with all pixels set to 1.
The operator $| \cdot |$ in this line  converts a vector-valued
image in $V$ to a magnitude image in $I$, and the $\max(\lambda \mathbf{1}_I, \cdot)$
operation thresholds the lower bound of the magnitude
image to $\lambda$ pixel-wise. Operationally, the division is performed by dividing the
spatial-vector at each pixel of the numerator by the scalar in the corresponding pixel of
the denominator.
Another potential source of confusion is computing the magnitude $\|(A,\nabla)\|_2$.
The power method for doing this is covered explicitly in Appendix \ref{sec:power}.
If it is desired to enforce the positivity constraint, the indicator $\delta_P(u)$ can be added
to the primal objective, and the effect of this is indicator is the same as for Listing \ref{alglspos};
namely the right hand side of line \ref{ulstvupdate} goes inside the $pos(\cdot)$ operator.

\subsubsection{Alternate data divergences}

For a number of reasons motivated by the physical model of imaging systems, it may be of use
to formulate optimization problems for CT image reconstruction with alternate data-error terms.
A natural extension of the quadratic data divergence is to include a diagonal weighting matrix.
The corresponding CP algorithm instance can be easily derived following the steps mentioned above.
As pointed out above, the CP method is not limited to quadratic objective functions and
other important convex functions can be used. We derive, here, three additional CP algorithm
instances.
For alternate data divergences we consider the oft-used KL divergence,
and one not so commonly used, $\ell_1$ data-error norm.
For the following, we need only
analyze the function $F_1$, as everything else remains the same as for the $\ell_2^2$-TV
objective in Eq. (\ref{LSTVprimal}).

\paragraph*{TV plus KL data divergence}
One data divergence of particular interest for tomographic image reconstruction is KL.
Objectives based on KL are what is being optimized in the various forms
of MLEM, and it is used often when data noise is a significant physical factor
and the data are modeled as being drawn from a multivariate
Poisson probability distribution \cite{Barrett:FIS}. For the situation where the view-sampling
is also sparse, it might be of interest to 
combine a KL data error term with the TV semi-norm in the following primal optimization:
\begin{linenomath}
\begin{equation}
\label{KLTVprimal}
\min_u \; \left\{ 
\sum_i \left[ Au - g+ g \ln g - g \ln(pos(Au))\right]_i + \delta_P(Au) + \lambda \left\| (|\nabla u|) \right\|_1 \right\},
\end{equation}
\end{linenomath}
where $\sum_i [\cdot]_i$ performs summation over all components of the vector argument.
This example proceeds as above except that the $F_1$ function is different: 
\begin{linenomath}
\begin{equation}
\label{KLF}
F_1(y) = \sum_i \left[ y -g + g \ln g - g \ln(pos(y))\right]_i + \delta_P(y)
\end{equation}
\end{linenomath}
where $y \in D$, and the function $\ln$ operates on the components of its argument.
Use of the KL data divergence makes sense only with positive linear systems $A$ and
non-negative pixel values $u$ and data $g$.
However, by defining the function over the whole space and using an indicator
function to restrict the domain \cite{rockafellar1970convex}, a wide variety of optimization
problems can be treated in a uniform manner.
Accordingly, $\delta_P$ is introduced into the $F_1$
objective and the $pos$ operator is used just so that this objective is defined
 in the real numbers.
% but
%the strategy for convex analysis promoted in Ref. \cite{rockafellar1970convex} is
%to define functions over all space and use indicator functions to restrict the domain.
%Accordingly, the $\delta_P$ is introduced into the $F_1$ objective and the $pos$ operator
%is used just so that this objective is defined in the real numbers.
%The fact that functions, like $F_1$, are defined for all $y$ allows for a
%wide variety of optimization problems to be treated in a
%uniform manner.
The derivation of $F_1^*$, though mechanical, is a little bit too long to be included
here. We simply state the resulting conjugate function:
\begin{linenomath}
\begin{equation}
\label{KLFconj}
F_1^*(p) =  \sum_i \left[ - g \ln pos( \mathbf{1}_D - p) \right]_i
+\delta_P(\mathbf{1}_D - p) .
\end{equation}
\end{linenomath}
The resulting dual problem to Eq. (\ref{KLTVprimal}) is thus:
\begin{linenomath}
\begin{equation}
\label{KLTVdual}
\max_{p,q} \; \left\{ \sum_i \left[ g \ln pos( \mathbf{1}_D - p) \right]_i
-\delta_P(\mathbf{1}_D - p) - \delta_{Box(\lambda)}(|q|) - \delta_{\mathbf{0}_I}(div\, q- A^T p) \right\}.
\end{equation}
\end{linenomath}

To form the algorithm instance, we need the proximal
mapping $prox_\sigma[F_1^*](y)$
\begin{linenomath}
\begin{equation}
\label{KLprox}
prox_\sigma[F_1^*](y) = \frac{1}{2}
\left( \mathbf{1}_D + y - \sqrt{ (y - \mathbf{1}_D)^2 + 4 \sigma g} \right).
\end{equation}
\end{linenomath}
An interesting point in the derivation, shown partially in Appendix \ref{sec:proximals},
of $prox_\sigma[F_1^*](y)$ is that the quadratic equation is needed, and the support function in
$F_1^*(p)$ is used to select the correct
(in this case negative) root of the discriminant in the quadratic formula.
With the new function $F_1$, its conjugate, and the conjugate's proximal mapping, we
can write down the CP algorithm instance.
Listing \ref{algkltv} gives the CP algorithm instance minimizing a KL plus
TV semi-norm objective.
\begin{algorithm}
\caption{Pseudocode for $N$-steps of the KL-TV CP algorithm instance.}
\label{algkltv}
\begin{algorithmic}[1]
\STATE $L \gets \|(A,\nabla)\|_2; \; \tau \gets 1/L; \; \sigma \gets 1/L; \; \theta \gets 1; \; n \gets 0$
\STATE initialize $u_0$, $p_0$, and $q_0$ to zero values
\STATE $\bar{u}_0 \gets u_0$
\REPEAT
\STATE $p_{n+1} \gets \frac{1}{2}
\left( \mathbf{1}_D + p_n+\sigma A\bar{u}_n - \sqrt{ ( p_n+\sigma A\bar{u}_n- \mathbf{1}_D)^2 + 4 \sigma g} \right)$
\label{pkltvupdate}
\STATE $q_{n+1} \gets \lambda (q_n + \sigma \nabla \bar{u}_n )/
\max(\lambda \mathbf{1}_I,|q_n + \sigma \nabla \bar{u}_n |)$ \label{qkltvupdate}
\STATE $u_{n+1} \gets  u_n - \tau A^T p_{n+1} + \tau div \, q_{n+1}$ \label{ukltvupdate}
\STATE $\bar{u}_{n+1} \gets u_{n+1} + \theta(u_{n+1} - u_n)$
\STATE $n \gets n+1$
\UNTIL{$n \ge N$}
\end{algorithmic}
\end{algorithm}
The difference between this algorithm instance and the previous $\ell_2^2$-TV
case comes only at the update at line \ref{pkltvupdate}.
This algorithm instance has the interesting property that the intermediate
image estimates $u_n$ can have negative values even though
the converged solution will be non-negative.
If it is desirable to have the intermediate image estimates be non-negative, the non-negativity
constraint can be easily introduced by adding the indicator $\delta_P(u)$ to the primal
objective, resulting in the addition of the $pos(\cdot)$ operator at
line \ref{ukltvupdate} as was shown in Listing \ref{alglspos}.

\paragraph*{TV plus $\ell_1$ data-error norm}
The combination of TV semi-norm regularization and $\ell_1$ data-error norm has been proposed
for image denoising and it has some interesting properties for that purpose \cite{chan2005aspects}.
This objective is also presented in Ref. \cite{chambolle2011first}. For tomography, this combination
may be of interest because the $\ell_1$ data-error term is an example of a robust fit to
the data. The idea of robust approximation is to weakly penalize data that are outliers \cite{boyd2004convex}.
Fitting with the commonly used quadratic error function, clearly puts heavy weight on outlying
measurements which in some situations can lead to streak artifacts in the images.  In particular,
for tomographic image reconstruction with a pixel basis, discretization error and metal objects can lead
to highly non-uniform error in the data model. Use of the $\ell_1$ data-error term may allow
for large errors for measurements along the tangent rays to internal structures, where discretization
can have a large effect. The $\ell_1$ data-error term also puts greater emphasis on fitting
the data that lie close to the model. The primal problem of interest is:
\begin{linenomath}
\begin{equation}
\label{L1TVprimal}
\min_u \; \left\{ \|Au - g \|_1  +\lambda \left\| (|\nabla u|) \right\|_1 \right\}.
\end{equation}
\end{linenomath}
For this objective, the function $F_1$ is:
\begin{linenomath}
\begin{equation}
\label{L1F}
F_1(y) = \| y - g \|_1.
\end{equation}
\end{linenomath}
Computing the convex conjugate $F_1^*$ yields
\begin{linenomath}
\begin{equation}
\label{L1Fconj}
F_1^*(p) = \delta_{Box(1)}(p) + \langle p, g \rangle_D,
\end{equation}
\end{linenomath}
and the resulting dual problem is:
\begin{linenomath}
\begin{equation}
\label{L1TVdual}
\max_{p,q}\;  \left\{ -\delta_{Box(1)}(p) - \langle p, g \rangle_D 
- \delta_{Box(\lambda)}(|q|) - \delta_{\mathbf{0}_I}(div \, q- A^T p) \right\}.
\end{equation}
\end{linenomath}

The proximal mapping necessary for completing the algorithm instance is
\begin{linenomath}
\begin{equation}
\label{L1prox}
prox_\sigma[F_1^*](y) = \frac{y-g\sigma}{\max(\mathbf{1}_D, |y-g \sigma|)},
\end{equation}
\end{linenomath}
where $\mathbf{1}_D$ is a data array with each component set to one and the $\max$ operation is
performed component-wise.
The corresponding pseudo-code for minimizing Eq. (\ref{L1TVprimal}) is given
in Listing \ref{algl1tv}, where the only difference between this code and the
previous two occurs at line \ref{pl1tvupdate}.
\begin{algorithm}
\caption{Pseudocode for $N$-steps of the $\ell_1$-TV CP algorithm instance.}
\label{algl1tv}
\begin{algorithmic}[1]
\STATE $L \gets \|(A,\nabla)\|_2; \; \tau \gets 1/L; \; \sigma \gets 1/L; \; \theta \gets 1; \; n \gets 0$
\STATE initialize $u_0$, $p_0$, and $q_0$ to zero values
\STATE $\bar{u}_0 \gets u_0$
\REPEAT
\STATE $p_{n+1} \gets (p_n+\sigma( A\bar{u}_n -g))/\max( \mathbf{1}_D, |p_n+\sigma( A\bar{u}_n -g)|)$
\label{pl1tvupdate}
\STATE $q_{n+1} \gets \lambda (q_n + \sigma \nabla \bar{u}_n )/ \max(\lambda \mathbf{1}_I,|q_n + \sigma \nabla \bar{u}_n |)$ \label{ql1tvupdate}
\STATE $u_{n+1} \gets  u_n - \tau A^T p_{n+1} + \tau div \, q_{n+1}$ \label{ul1tvupdate}
\STATE $\bar{u}_{n+1} \gets u_{n+1} + \theta(u_{n+1} - u_n)$
\STATE $n \gets n+1$
\UNTIL{$n \ge N$}
\end{algorithmic}
\end{algorithm}
The ability to deal with non-smooth objectives uncomplicates this particular
problem substantially. If smoothness were required, there would have to be
smoothing parameters on both the $\ell_1$ and TV terms, adding two more
parameters than necessary to a study of the image properties as a function of
the optimization-problem parameters.

\subsubsection{Constrained, TV-minimization}
The previous three optimization problems combine a data fidelity term
with a TV-penalty, and the balance of the two terms is controlled by
the parameter $\lambda$.  An inconvenience of such optimization
problems is that it is difficult to physically interpret $\lambda$.
Focusing on combining an $\ell_2$ data-error norm with TV, reformulating
Eq. (\ref{LSTVprimal}) as a constrained, TV-minimization leads
to the following primal problem:
\begin{linenomath}
\begin{equation}
\label{constrainedTVprimal}
\min_u \; \left\{ \|(|\nabla u|)\|_1 +\delta_{Ball(\epsilon)}(Au-g) \right\},
\end{equation}
\end{linenomath}
where $\delta_{Ball(\epsilon)}(Au-g)$ is zero for $\|Au-g\|_2 \le \epsilon$.
When $\epsilon>0$, this problem is equivalent to the unconstrained optimization
Eq. (\ref{LSTVprimal}), see e.g. Ref. \cite{elad2010sparse},
in the sense that for each positive $\epsilon$ there is a corresponding $\lambda$
yielding the same solution.
For this constrained, TV-minimization, the function $F_1$ is
\begin{linenomath}
\begin{equation}
\label{ConstrainedTVF}
F_1(y) = \delta_{Ball(\epsilon)}(y -g).
\end{equation}
\end{linenomath}
The corresponding conjugate is
\begin{linenomath}
\begin{equation}
\label{ConstrainedTVFconj}
F_1^*(p) = \epsilon \| p \|_2 + \langle p,g \rangle_D,
\end{equation}
\end{linenomath}
leading to the dual problem:
\begin{linenomath}
\begin{equation}
\label{constrainedTVdual}
\max_{p,q} \left\{- \epsilon \| p \|_2 - \langle p,g \rangle_D - \delta_{Box(1)}(|q|)
- \delta_{\mathbf{0}_I}(div \, q- A^T p) \right\}.
\end{equation}
\end{linenomath}

Again,
for the algorithm instance we need the proximal mapping $prox_\sigma[F_1^*]$:
\begin{linenomath}
\begin{equation}
\label{ConstrainedTVprox}
prox_\sigma[F_1^*](y) = \max( \|y- \sigma g \|_2 - \sigma \epsilon, 0) \, (y- \sigma g).
\end{equation}
\end{linenomath}
The main points in deriving this proximal mapping are discussed in Appendix \ref{sec:proximals},
and it is an example where geometric/symmetry arguments play a large role.
Listing \ref{algConstrainedTV} shows the algorithm instance solving Eq. (\ref{constrainedTVprimal}),
where once again only line \ref{pConstrainedTVupdate} is modified.
\begin{algorithm}
\caption{Pseudocode for $N$-steps of the $\ell_2$-constrained, TV-minimization CP algorithm instance.}
\label{algConstrainedTV}
\begin{algorithmic}[1]
\STATE $L \gets \|(A,\nabla)\|_2; \; \tau \gets 1/L; \; \sigma \gets 1/L; \; \theta \gets 1; \; n \gets 0$
\STATE initialize $u_0$, $p_0$, and $q_0$ to zero values
\STATE $\bar{u}_0 \gets u_0$
\REPEAT
\STATE $p_{n+1} \gets \max( \| p_n+\sigma( A\bar{u}_n -g)\|_2 - \sigma \epsilon, 0 ) \, (p_n+\sigma( A\bar{u}_n -g))$
\label{pConstrainedTVupdate}
\STATE $q_{n+1} \gets (q_n + \sigma \nabla \bar{u}_n )/ \max(\mathbf{1}_I,|q_n + \sigma \nabla \bar{u}_n |)$ \label{qConstrainedTVupdate}
\STATE $u_{n+1} \gets  u_n - \tau A^T p_{n+1} + \tau div \, q_{n+1}$ \label{uConstrainedTVupdate}
\STATE $\bar{u}_{n+1} \gets u_{n+1} + \theta(u_{n+1} - u_n)$
\STATE $n \gets n+1$
\UNTIL{$n \ge N$}
\end{algorithmic}
\end{algorithm}
This algorithm instance essentially achieves the same goal as Listing \ref{alglstv}, the only difference
is that the parameter $\epsilon$  has an actual physical interpretation, being the data-error
bound.

\section{Demonstration of CP algorithm instances for tomographic image reconstruction}
\label{sec:simulations}

In the previous section, we have derived CP algorithm instances covering many optimization
problems of interest to CT image reconstruction. Not only are there the seven optimization
problems, but within each case the system model/matrix $A$, the data $g$, and optimization
problem parameters can vary.  For each of these, practically infinite number of optimization
problems, the corresponding CP algorithm instances are guaranteed to converge
\cite{chambolle2011first}. The purpose
of this results section is not to advocate one optimization problem over another; rather to
demonstrate the utility of the CP algorithm for optimization problem prototyping. 
For this purpose, we present example image reconstructions that could be performed
in a study 
for investigating the impact of matching the data-divergence with data noise model for 
image reconstruction in breast CT.

\subsection{Sparse-view experiments for image reconstruction from simulated CT data}

We briefly describe the significance of the experiments, but we point out that
the main goal here is to demonstrate the CP algorithm instances.
Much of the recent interest in employing the TV semi-norm in optimization problems for CT image
reconstruction has been generated by compressive sensing (CS).  CS seeks to relate sampling
conditions on a sensing device with sparsity in the object being scanned. So far, mathematical
results have been limited to various types of random sampling \cite{candes2008introduction}.  System matrices
such as those representing CT projection fall outside of the scope of mathematical results
for CS \cite{SidkyPC:10}. As a result, the only current option for investigating CS in CT is through numerical
experiments with computer phantoms.

%Recently, we have been interested in characterizing CS sampling conditions for breast CT \cite{jakob:2011}.
%Preliminarily, it seems that constrained, TV minimization is fairly robust with a circular scanning
%configuration with equi-angular spacing.  Numerical experiments seem to show that the necessary number
%of samples for accurate recovery under ideal conditions is around 2.5 times the object sparsity, which
%is fairly close to the theoretical limit of twice the object sparsity \cite{candes2006robust}.
%These experimental results paint a far more optimistic picture about exploiting sparsity than
%our attempts at proving reconstructability through the restricted isometry property \cite{SidkyPC:10}.
%The controlled study
%of Ref. \cite{jakob:2011} together with multiple successes with actual CT scanner data,
%indicate that sparsity-exploiting image reconstruction for CT may perform better in practice than in theory
%-- a rare circumstance.

A next logical step for bridging theoretical results for CS
to actual application is to consider physical
factors in the data model.  One such factor is a noise model, which can be quite important for low-dose
CT applications such as breast CT.  While much work has been performed on iterative image reconstruction
with various noise models under conditions of full sampling, little is known about the impact of noise
on sparse-view image reconstruction.  In the following limited study, we set up a breast CT simulation
to investigate the impact of correct modeling of data noise with the purpose of
demonstrating that the CP algorithm instances can be applied to the CT system.

\subsection{Sparse-view reconstruction with a Poisson noise model}

\begin{figure}[!h]
\begin{minipage}[b]{\linewidth}
\centering
\centerline{\includegraphics[width=\linewidth]{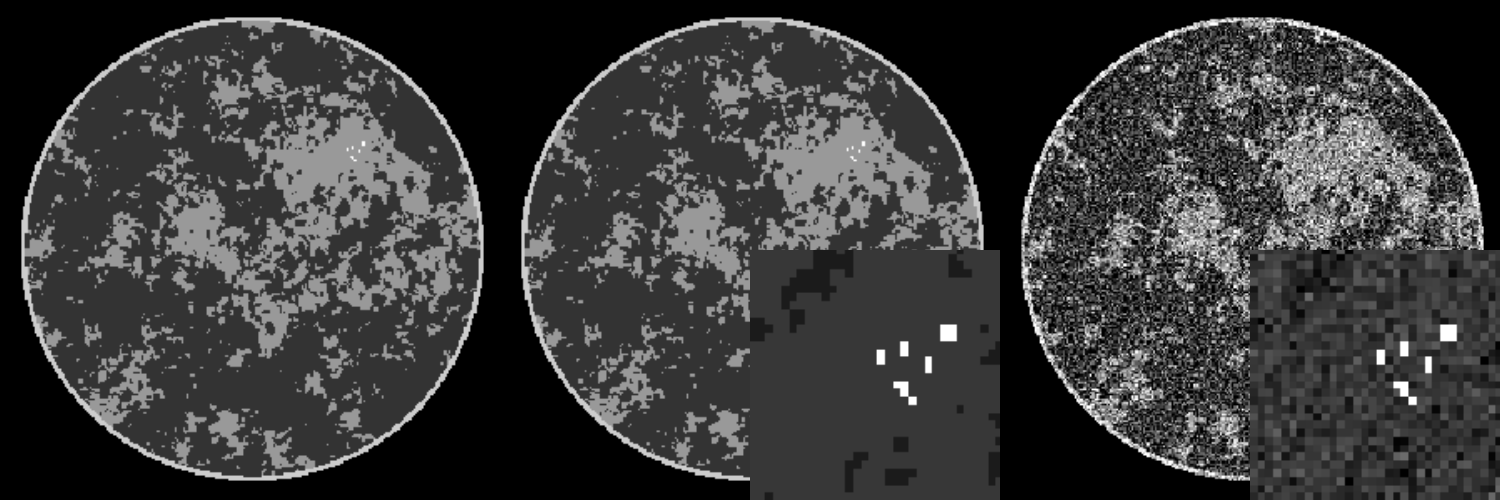}}
% \centerline{\includegraphics[width=8cm,clip=TRUE]{figs/viewscond.eps}}
\end{minipage}
\caption{Breast phantom for CT and FBP reconstructed image for a 512-view
data set with Poisson distributed noise.
%simulating a total detected photon flux of $10^{10}$ photons.
Left is the phantom in the gray scale window [0.95,1.15]; middle is the same
phantom with a blow-up on the micro-calcification ROI displayed in the gray
scale window [0.9,1.8]; and right is the FBP image reconstructed from the noisy data.
The middle panel is the reference for all image reconstruction algorithm results.
The FBP image is shown only to provide a sense of the noise level.
\label{fig:phantom}}
\end{figure}
For the following study, we employ a digital 256 $\times$ 256 breast phantom, described in
Ref. \cite{joergensen2011toward,reiser2010task},
and used in our previous study on investigating sufficient sampling conditions
for TV-based CT image reconstruction \cite{jakob:2011}.  The phantom models for tissue types:
the background fat tissue is a assigned a value of 1.0,
the modeled fibro-glandular tissue takes a value of 1.1, the outer skin layer is set to 1.15,
and the micro-calcifications are assigned values in the range [1.8,2.3].

For the present case, we focus on circular, fan-beam
scanning with 60 projections equally distributed over a full 360$^\circ$ angular range.
The simulated radius of the X-ray source trajectory is 40cm with a source-detector distance
of 80cm.
The detector sampling consists of 512 bins of size 200 microns.
The system matrix for the X-ray projection is computed by the line-intersection method
where the matrix elements of $A$ are determined by the length of traversal in each
image pixel of each source/detector-bin ray. 
For this phantom under ideal conditions, we
have found that accurate recovery is possible with constrained, TV-minimization with
as few as 50 projections. In the present study, we add Poisson noise to the data model
at a level consistent with what might be expect in a typical breast CT scan. The Poisson
noise model is chosen in order to investigate the impact of matching the data-error term
to the noise model.
For reference, the phantom is shown in Fig. \ref{fig:phantom}. To have a sense of the noise level,
a standard fan-beam filtered back-projection reconstruction is shown alongside the phantom for
simulated Poisson noise.% with $10^{10}$ counts for a 512-view data set.

For this noise model, the maximum likelihood method prescribes minimizing the KL data divergence
between the available and estimated data. To gauge the importance of selecting a maximum likelihood
image, we compare the results from two optimization problems:
a KL data divergence plus a TV-penalty, Eq. (\ref{KLTVprimal}) above; and a least-squares data error norm
plus a TV-penalty, Eq. (\ref{LSTVprimal}) above. With the CP framework, these two optimization problems
can be easily prototyped: the solutions to both problems can be obtained without worrying about smoothing
the TV semi-norm, setting algorithm parameters, or proving convergence.

\begin{figure}[!h]
\begin{minipage}[b]{\linewidth}
\centering
\centerline{\includegraphics[width=\linewidth]{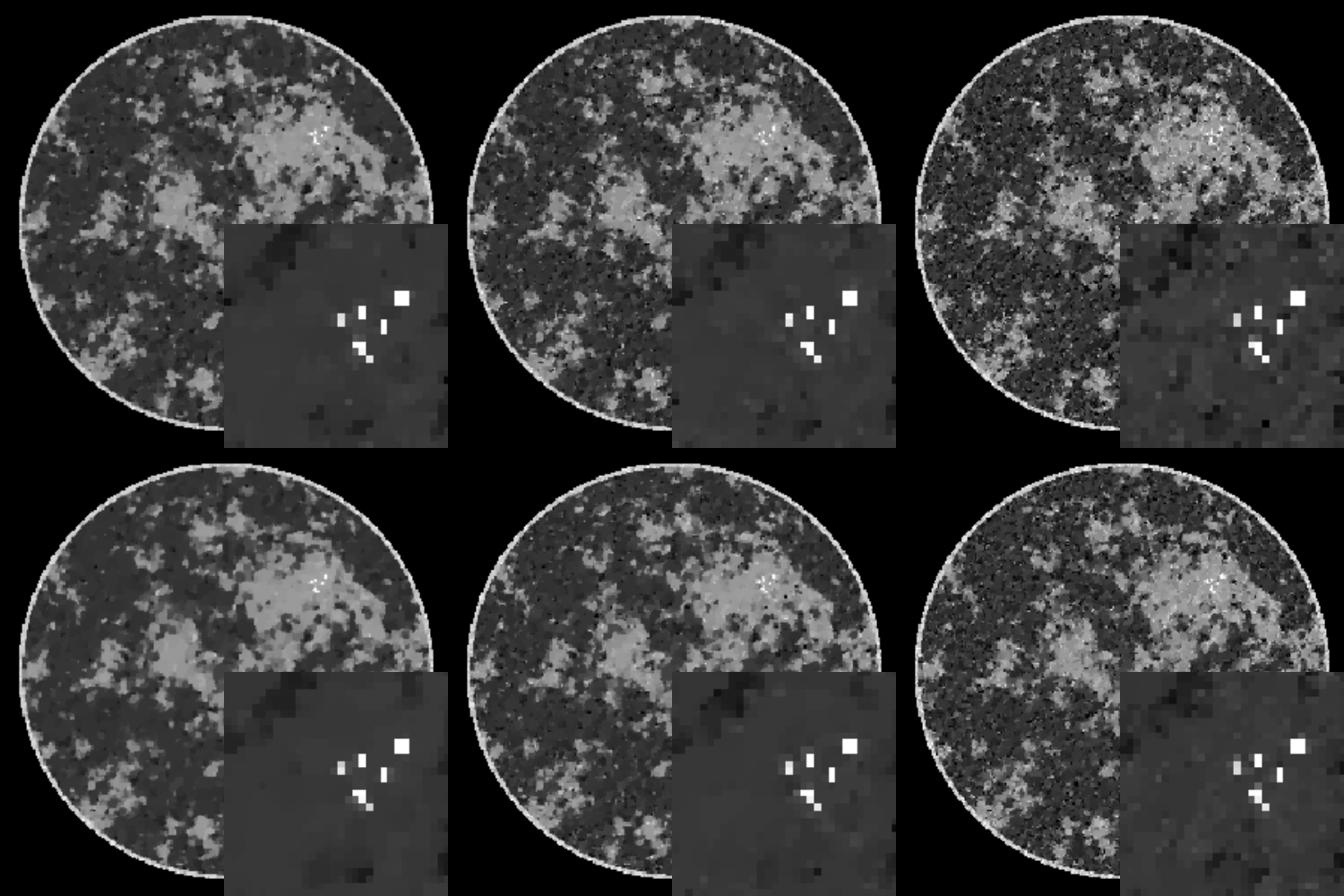}}
% \centerline{\includegraphics[width=8cm,clip=TRUE]{figs/viewscond.eps}}
\end{minipage}
\caption{
Images reconstructed from 60-view projection data with a Poisson distributed noise model.
The top row of images result from minimizing
the $\ell^2_2$-TV objective in Eq. (\ref{LSTVprimal}) for $\lambda= 1 \times 10^{-4}$, $5 \times 10^{-5}$,
and $2 \times 10^{-5}$, going from left to right. The bottom row of images result from minimizing
the KL-TV objective in Eq. (\ref{KLTVprimal}) for the same values of $\lambda$.
Note that $\lambda$ does not necessarily have the same impact on each of these
optimization problems. Nevertheless, we see similar trends for the chosen
values of $\lambda$.
\label{fig:images}}
\end{figure}
For the phantom and data conditions, described above, the images for different values of the TV-penalty
parameter $\lambda$ are shown in Fig. \ref{fig:images}. An ROI of the micro-calcification cluster is also shown.
The overall and ROI images give an impression of two different visual tasks important for breast imaging:
discerning the fibro-glandular tissue morphology and detection/classification of micro-calcifications.
The images show some difference between the two optimization problems;
most notably there is a perceptible reduction in noise in the ROIs from the
KL-TV images. A firm conclusion, however, awaits a more complete study with multiple noise realizations.

The most critical feature of the CP algorithm that we wish to promote is the rapid prototyping of
a convex optimization problem for CT image reconstruction. The above study is aimed at a combination
of using a data divergence based on maximum likelihood estimation with a TV-penalty, which takes
advantage of sparsity in the gradient magnitude of the underlying object.  The CP framework facilitates
the use of many other convex optimization problems, particularly those based on some form of sparsity,
which often entail some form of the non-smooth $\ell_1$-norm. For example, in Ref. \cite{SidkyPC:10}
we have found it useful for sparse-view X-ray phase-contrast imaging to perform image reconstruction with
a combination of a least-squares data fidelity term, an $\ell_1$-penalty promoting
object sparseness, and an image TV constraint to further reduce streak artifacts from angular
under-sampling. Under the CP framework, prototyping various combinations of these terms as
constrained or unconstrained optimization problems becomes possible and the corresponding
derivation of CP algorithm instances follows from the steps described in Sec. \ref{sec:prototyping}.
Alternative, convex data fidelity terms and image constraints motivated by
various physical models may also be prototyped.

As a practical matter, though, it is important to have some sense of the convergence
of the CP algorithm instances.
To this end, we take an in depth look at individual runs 
for the KL-TV algorithm instance for CT image reconstruction.

\subsection{Iteration dependence of the CP algorithm}

\begin{figure}[!h]
\begin{minipage}[b]{0.49\linewidth}
\centering
\centerline{\includegraphics[width=0.95\linewidth]{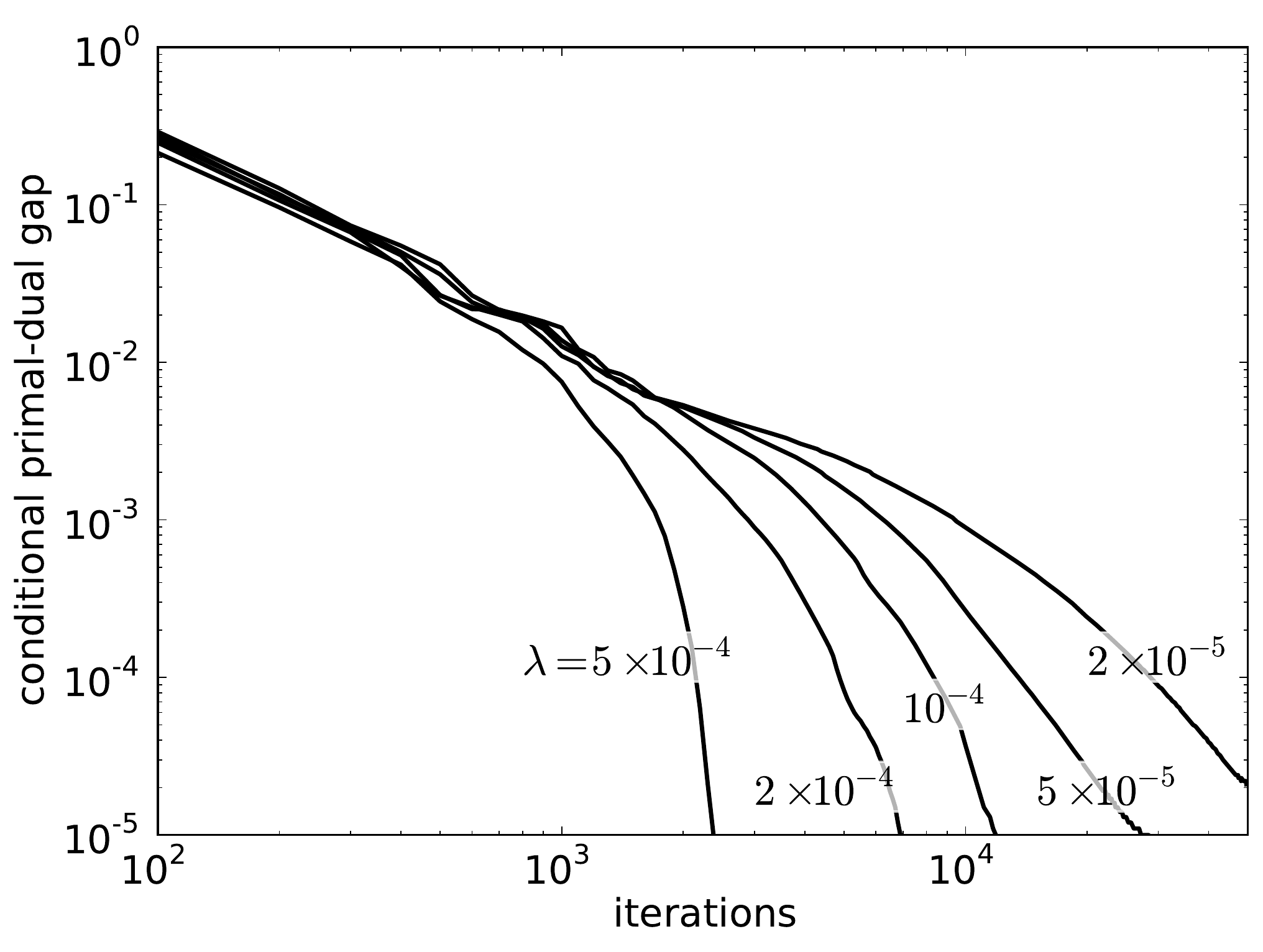}}
% \centerline{\includegraphics[width=8cm,clip=TRUE]{figs/viewscond.eps}}
\end{minipage}
\begin{minipage}[b]{0.49\linewidth}
\centering
\centerline{\includegraphics[width=0.95\linewidth]{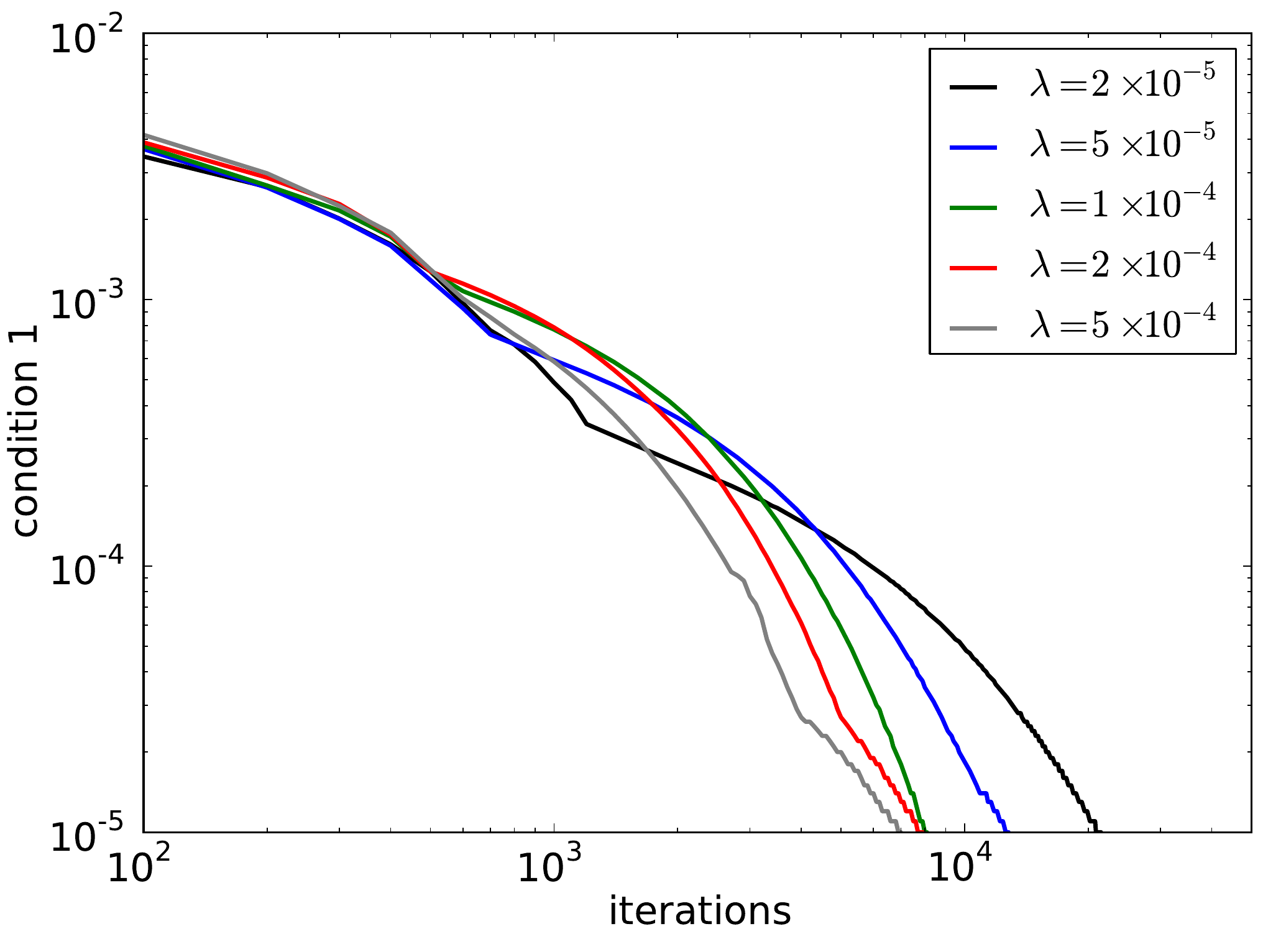}}
% \centerline{\includegraphics[width=8cm,clip=TRUE]{figs/viewscond.eps}}
\end{minipage}
\caption{
(Left) Convergence of the partial primal-dual gap for the CP algorithm instance solving 
Eq. (\ref{KLTVprimal}) for different values of $\lambda$.
(Right) Plot indicating agreement with
condition 1: $\| div\, q- A^T p \|_\infty$, the magnitude of the 
largest component of the argument of the last indicator function of Eq. (\ref{KLTVdual}).
Collecting all the indicator functions of the primal, Eq. (\ref{KLTVprimal}), and
dual, Eq. (\ref{KLTVdual}), KL-TV optimization problems, we have four conditions
to check in addition to the conditional primal-dual gap: (1) $div\, q- A^T p=\mathbf{0}_I$,
(2) $Au \geq\mathbf{0}_D$, (3) $p\leq \mathbf{1}_D$, and (4) $|q|<\lambda$.
The agreement with condition 1 is illustrated in the plot; agreement with condition 2 
has a similar dependence; condition 3 is satisfied early on in the iteration; and
condition 4 is automatically enforced by the CP algorithm instance for KL-TV.
Because the curves
are bunched together in the condition 1 plot, they are differentiated in color.
\label{fig:pdgap}}
\end{figure}

Through the methods described above, many
useful algorithm instances can be derived for CT image reconstruction.
It is obviously important that the
resulting algorithm instance reaches the solution of the prescribed optimization problem.
To illustrate the convergence of a resulting algorithm instance we focus on the
TV-penalized KL data divergence, Eq. (\ref{KLTVprimal}), and plot
the conditional
primal-dual gap for the different runs with varying $\lambda$ in Fig. \ref{fig:pdgap}.
Included in this figure is a plot indicating the convergence to agreement with
the most
challenging condition set by the indicator functions in Eq. (\ref{KLTVdual}).
For the present results we terminated the iteration at a conditional primal-dual gap
of $10^{-5}$,
which appears to happen on the scale of thousands of iterations with
smaller $\lambda$ requiring more iterations. Interestingly, a simple pre-conditioned
form of the CP algorithm was proposed in Ref. \cite{Pock2011}, which appears to
perform efficiently for small $\lambda$. The pre-conditioned CP algorithm instance for this
problem is reported in Appendix \ref{sec:preconditionedCP}.

%  To place the convergence of the conditional primal-dual
%gap in context, we also plot the convergence of the RMSE, ROI-RMSE, and gradient-RMSE for the
%smallest $\lambda$.
%The conditional primal-dual gap provides a powerful test of convergence,
%because it does not require the objective functions to be smooth.

\section{Discussion}

This article has presented the application of the Chambolle-Pock algorithm to prototyping
of optimization problems for CT image reconstruction.  The algorithm covers many optimization
problems of interest allowing for non-smooth functions. It also comes with solid convergence
criteria to check the image estimates.

The use of the CP algorithm we are promoting here is for prototyping;
namely, when the image reconstruction algorithm development is at the early
stage of determining important factors in formulating the optimization problem.
As an example, we illustrated a scenario for sparse-view breast CT considering two
different data-error terms.  In this stage of development it is helpful to
not have to bother with algorithm parameters, and questions of whether or not
the algorithm will converge.  After the final optimization problem is determined,
then the focus shifts from prototyping to efficiency.

Optimization problem prototyping for CT image reconstruction does have its 
limitations.
For example, in the breast CT simulation presented above, a more complete
conclusion requires reconstruction from multiple realizations of the data under the Poisson
noise model. Additional important dimensions of the study are generation of an ensemble
of breast phantoms and considering alternate image representations/projector models.
Considering the size of CT image reconstruction systems and huge parameter space
of possible optimization problems, it is not yet realistic to completely characterize
a particular CT system. But at least we are assured of solving isolated setups
and it is conceivable to perform a study along one aspect of the system,
i.e. consider multiple realizations of the
random data model.  Given the current state of affairs for optimization-based
image reconstruction, it is crucial that simulations be as realistic as possible.
There is great need for realistic phantoms, and data simulation software.

%In considering efficiency there is a wide array of algorithm literature
%on first-order methods, including also potential
%efficiency gains with the Chambolle-Pock algorithm itself.
We point out that it is likely at least within the immediate future, that optimization-based
image reconstruction will have to operate at severely truncated iteration numbers.
Current clinical applications of iterative image reconstruction often operate in the
range of one to ten iterations, which is likely far too few for claiming that
the image estimate is an accurate solution to the designed optimization problem.
But at least the ability to prototype an optimization problem can potentially
simplify the design phase by separating optimization parameters from algorithm
parameters.

\section*{Acknowledgment}
This work is part of the project CSI: Computational Science
in Imaging, supported by grant 274-07-0065 from the Danish
Research Council for Technology and Production Sciences.
This work was supported in part by NIH R01 grants CA158446, CA120540, and EB000225.
The contents of this article are
solely the responsibility of the authors and do not necessarily
represent the official views of the National Institutes of Health.

\appendices

\section{Computing the norm of $K$}
\label{sec:power}

The matrix norm used for the parameter $L$ in the CP algorithm instances is the largest
singular value of $K$.  This singular value can be obtained by the standard
power method specified in Listing \ref{power}.
\begin{algorithm}
\caption{Pseudocode for $N$-steps of the generic power method.
The scalar $s$ tends to $\|K\|_2$ as $N$ increases.}
\label{power}
\begin{algorithmic}[1]
\STATE initialize $x_0 \in I$ to a non-zero image
\STATE $n \gets 0$
\REPEAT
\STATE $x_{n+1} \gets K^T K x_n$
\STATE $x_{n+1} \gets x_{n+1}/\|x_{n+1}\|_2$
\STATE $s \gets \| K x_{n+1} \|_2$
\STATE $n \gets n+1$
\UNTIL{$n \ge N$}
\end{algorithmic}
\end{algorithm}
When $K$ represents the discrete X-ray transform, our experience has
been that the power method converges to numerical precision in twenty
iterations or less. In implementing the CP algorithm instance for TV-penalized
minimization, the norm of the combined linear transform $\|(A,\nabla)\|_2$
is needed. For this case, the program is the same as Listing~\ref{power}
where $K^T Kx_n$ becomes $A^TAx_n -div \nabla x_n$; recall that $-div = \nabla^T$.
Furthermore, to obtain $s$, the explicit computation is
$s=\sqrt{\|A x_{n+1}\|^2_2 + \|\nabla x_{n+1}\|^2_2}$.

\section{The convex conjugate of certain indicator functions of interest
illustrated in one-dimension}
\label{sec:legendre}

This appendix covers the convex conjugate of a couple of indicator functions in
one dimension, serving to illustrate how geometry plays a role in the computation
and to provide a mental picture on the conjugate of
higher dimensional indicator functions.

\begin{figure}[!h]
\begin{minipage}[b]{\linewidth}
\centering
\centerline{\includegraphics[width=0.5\linewidth]{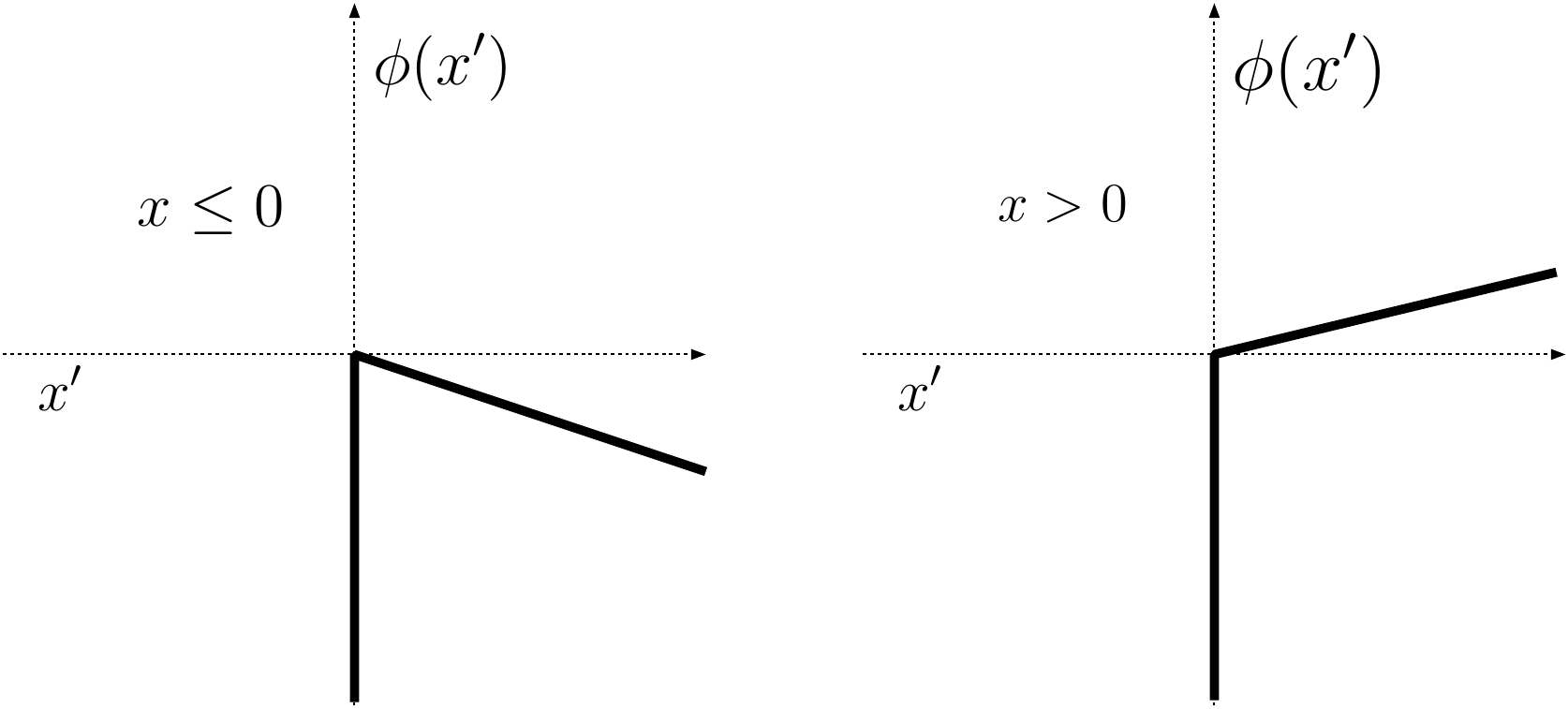}}
% \centerline{\includegraphics[width=8cm,clip=TRUE]{figs/viewscond.eps}}
\end{minipage}
\caption{
Illustration of the objective function, labeled $\phi (x^\prime)$, in the maximization
described by Eq. (\ref{leg_pos}). Shown are the two cases discussed in the text.
\label{fig:deltaP}}
\end{figure}
Consider first the indicator $\delta_P(x)$, which is zero for $x \ge 0$.
The conjugate of this indicator is computed from:
\begin{linenomath}
\begin{equation}
\label{leg_pos}
\delta^*_P(x) = \max_{x^\prime} \; \phi (x^\prime) =  \max_{x^\prime}\;
 \left\{x^\prime x - \delta_P(x^\prime) \right\}.
\end{equation}
\end{linenomath}
To perform this maximization, we analyze the cases, $x \le 0$ and $x > 0$, separately.
As a visual aid, we plot the objective for these two cases in Fig. \ref{fig:deltaP}.
From this figure it is clear that when $x \le 0$, the objective's maximum is attained
at $x^\prime=0$ and this maximum value is 0 (note that this is true even for $x=0$).
When $x>0$, the objective can increase without bound as $x^\prime$ tends to $\infty$,
resulting in a maximum value of $\infty$. Putting these two cases together yields:
\begin{linenomath}
\begin{equation}
\notag
\delta^*_P(x) = \delta_P(-x).
\end{equation}
\end{linenomath}
Generalizing this argument to multi-dimensional $x$, yields Eq. (\ref{legpos}).

\begin{figure}[!h]
\begin{minipage}[b]{\linewidth}
\centering
\centerline{\includegraphics[width=0.5\linewidth]{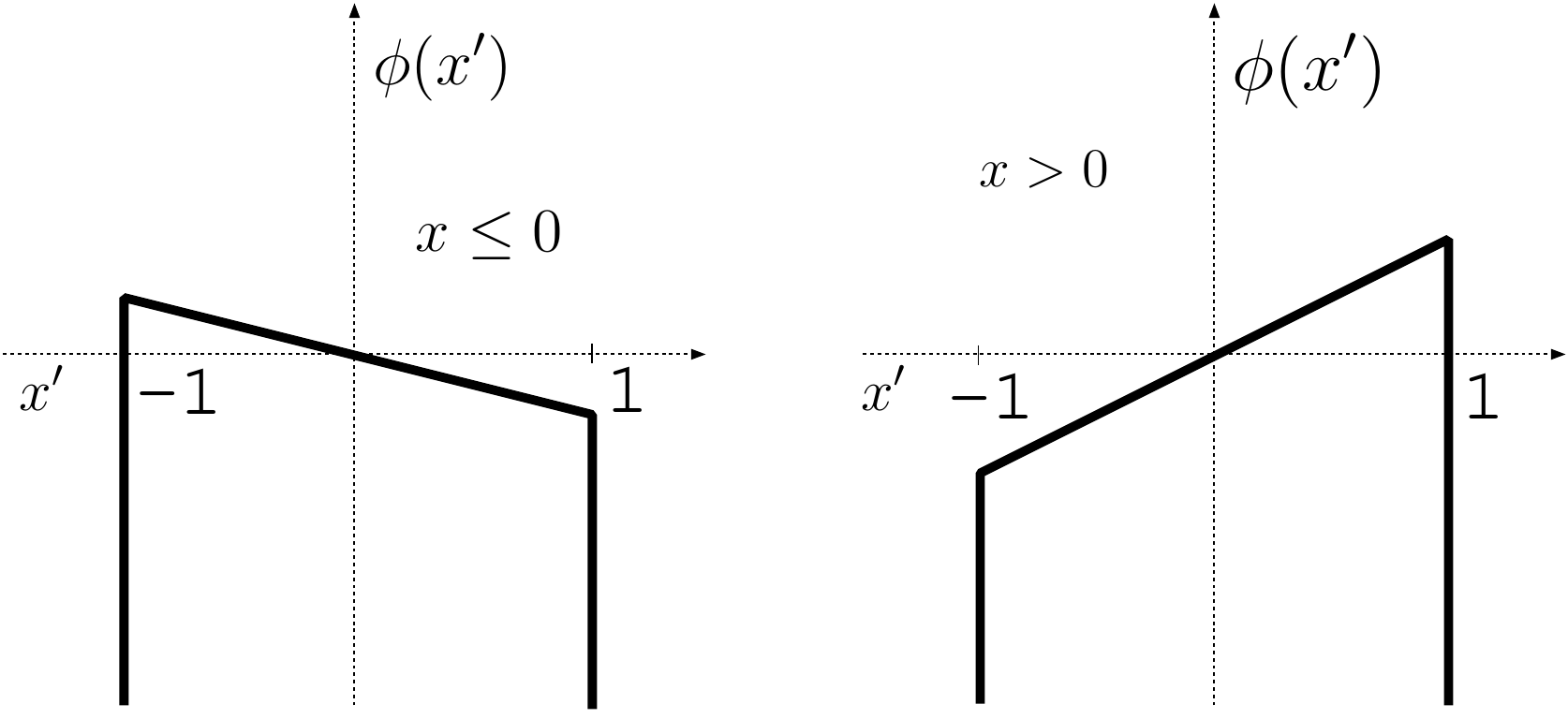}}
% \centerline{\includegraphics[width=8cm,clip=TRUE]{figs/viewscond.eps}}
\end{minipage}
\caption{
Illustration of the objective function, labeled $\phi (x^\prime)$, in the maximization
described by Eq. (\ref{leg_box}). Shown are the two cases discussed in the text.
\label{fig:deltaBox}}
\end{figure}
Next we consider $\delta_{Box(1)}(x)$, which in one dimension is the same as
$\delta_{Ball(1)}(x)$.  This functions is zero only for $-1 \le x \le 1$.
Its conjugate is computed from:
\begin{linenomath}
\begin{equation}
\label{leg_box}
\delta^*_{Box(1)}(x) =  \max_{x^\prime} \; \phi (x^\prime) =
 \max_{x^\prime} \; \left\{ x^\prime x - \delta_{Box(1)}(x^\prime) \right\}.
\end{equation}
\end{linenomath}
Again, we have two cases, $x \le 0$ and $x>0$, illustrated in Fig. \ref{fig:deltaBox}.
In the former case the maximum value of the objective is attained at $x^\prime = -1$,
and this maximum value is $-x$. In the latter case the maximum value is $x$, and
it is attained at $x^\prime = 1$. Hence, we have:
\begin{linenomath}
\begin{equation}
\notag
\delta^*_{Box(1)}(x) = |x|.
\end{equation}
\end{linenomath}
For multi-dimensional $x$, $\delta_{Box(1)}(x) \ne \delta_{Ball(1)}(x)$, and this is
also reflected in the conjugates:
\begin{linenomath}
\begin{align}
\notag
\delta^*_{Box(1)}(x) &= \|x\|_1, \\
\notag
\delta^*_{Ball(1)}(x) &= \|x\|_2.
\end{align}
\end{linenomath}

\begin{figure}[!h]
\begin{minipage}[b]{\linewidth}
\centering
\centerline{\includegraphics[width=0.25\linewidth,angle=270]{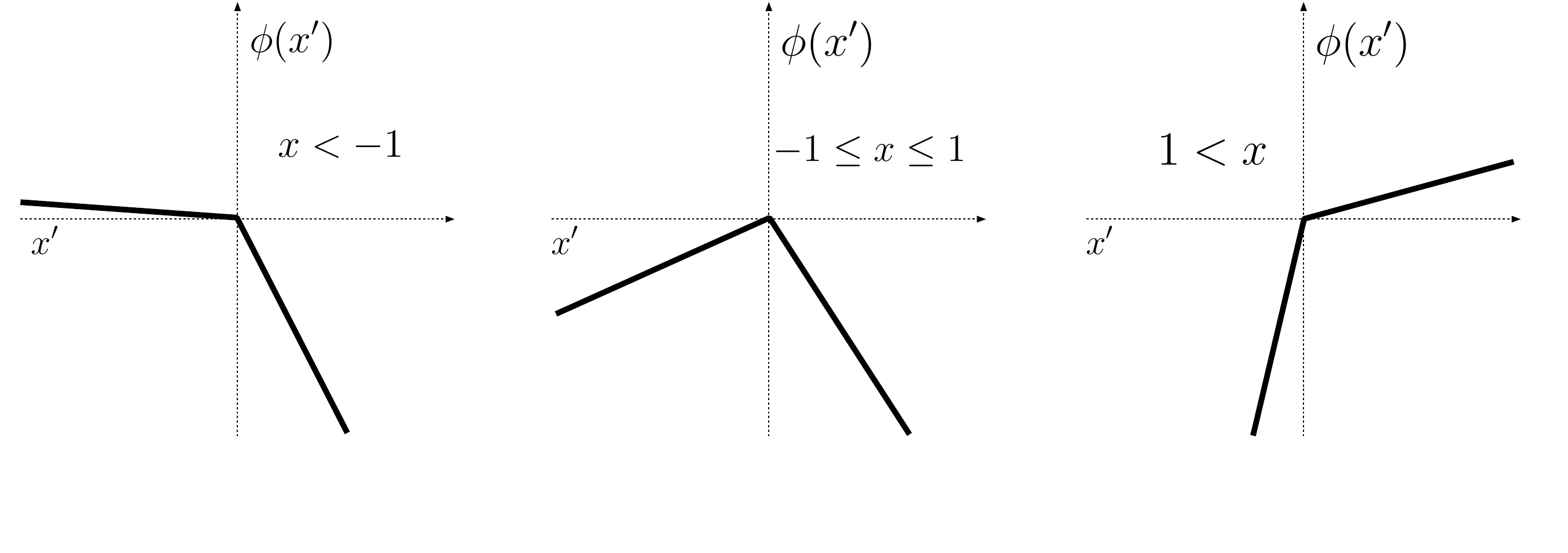}}
% \centerline{\includegraphics[width=8cm,clip=TRUE]{figs/viewscond.eps}}
\end{minipage}
\caption{
Illustration of the objective function, labeled $\phi (x^\prime)$, in the maximization
described by Eq. (\ref{leg_abs}). Shown are the three cases discussed in the text.
\label{fig:abs}}
\end{figure}
It is also interesting to verify that $\delta^{**}_{Box(1)}(x)$ is indeed $\delta_{Box(1)}(x)$
by showing, again in one dimension, that $|x|^* = \delta_{Box(1)}(x)$.
Illustrating this example helps in understanding the convex conjugate of
multi-dimensional $\ell_1$-based semi-norms. The relevant
conjugate is computed from:
\begin{linenomath}
\begin{equation}
\label{leg_abs}
|x|^* =  \max_{x^\prime} \; \phi (x^\prime) = \max_{x^\prime} \;
\left\{ x^\prime x - |x^\prime| \right\}.
\end{equation}
\end{linenomath}
Here, we need to analyze three cases: $x<-1$, $-1 \le x \le 1$, and $x>1$.
The corresponding sketch is in Fig. \ref{fig:abs}. The $-|x|$ term in the objective
makes an upside-down wedge, and the $x^\prime x$ term serves to tip this wedge.
In the second case, the wedge is tipped, but still opens up downward so that the objective
is maximized at $x^\prime=0$, attaining there the value of 0. In the first and third
cases, the wedge is tipped so much that part of it points upward and the objective
can increase without bound, attaining the value of $\infty$. Putting these cases together
does indeed yield:
\begin{linenomath}
\begin{equation}
\notag
|x|^* = \delta_{Box(1)}(x).
\end{equation}
\end{linenomath}
Similar reasoning is used to obtain Eq. (\ref{TVconj2}) from Eq. (\ref{TVconj}).

\section{Computation of important proximal mappings}
\label{sec:proximals}

This appendix fills in important steps in computing some of the proximal
mappings in the text, where it is necessary to use geometrical reasoning
in addition to setting the gradient of the objective to zero.

The conjugate of the TV semi-norm in Eq. (\ref{TVconj2}) leads to the following
proximal mapping computation:
\begin{linenomath}
\begin{equation}
\notag
prox_\sigma[F^*](z) = \argmin_{z^\prime}\; \left\{ \delta_{Box(\lambda)}
(|z^\prime|) + \frac{\| z - z^\prime \|_2^2}{2 \sigma} \right\},
\end{equation}
\end{linenomath}
where $z,z^\prime \in V$, and absolute value, $|\cdot|$, of a spatial-vector image $V$
yields an image, in $I$, of the spatial-vector-magnitude. The quadratic term
is minimized when $z=z^\prime$, but the indicator function excludes
this minimizer when $z \not\in Box(\lambda)$. To solve this problem,
we write the quadratic as a sum over pixels:
\begin{linenomath}
\begin{equation}
\notag
\frac{\| z - z^\prime \|_2^2}{2 \sigma} = \frac{ \sum_i |z_i-z_i^\prime|^2 }{2 \sigma},
\end{equation}
\end{linenomath}
where $i$ indexes the image pixels and each $z_i$ and $z_i^\prime$ is a spatial-vector.
The indicator function places an upper bound on the magnitude of
each spatial-vector $|z_i^\prime| \le \lambda$. The proximal mapping is
built pixel-by-pixel considering two cases: if $|z_i| \le \lambda$,
then $prox_\sigma[F^*](z)_i = z_i$; if $|z_i| > \lambda$, then $z_i^\prime$ is
chosen to be closest to $z_i$ while respecting $|z_i^\prime| \le \lambda$
which leads to a scaling of the magnitude of $z_i$ and
$prox_\sigma[F^*](z)_i = \lambda z_i/|z_i|$. Note that the constant
$\sigma$ does not enter into this calculation. Putting the cases and
components all together
yields the second part of the proximal mapping in Eq. (\ref{LSTVprox}).

For the KL-TV problem the proximal mapping for the data term is computed
from Eq. (\ref{KLFconj}):
\begin{linenomath}
\begin{equation}
\notag
prox_\sigma[F_1^*](p) = \argmin_{p^\prime}\;
\left\{
\frac{\| p - p^\prime \|_2^2}{2 \sigma}
-  \sum_i \left[ g \ln pos( \mathbf{1}_D - p^\prime)\right]_i
+\delta_P(\mathbf{1}_D - p^\prime) \right\} .
\end{equation}
\end{linenomath}
We note that the objective is a smooth function in the positive orthant
of $p^\prime \in D$. Accordingly, we differentiate the objective with
respect to $p^\prime$ ignoring the $pos(\cdot)$ and indicator functions,
keeping in mind that we have to check that the minimizer $p^\prime$ is
non-negative.
Performing the differentiation and setting to zero yields the following
quadratic equation:
\begin{linenomath}
\begin{equation}
\notag
{p^\prime}^2 - (\mathbf{1}_D+p) p^\prime +p - \sigma g  =0,
\end{equation}
\end{linenomath}
and substituting into the quadratic equation yields:
\begin{linenomath}
\begin{equation}
\notag
 p^\prime = \frac{1}{2} (\mathbf{1}_D + p  \pm \sqrt{ (\mathbf{1}_D-p)^2 + 4 \sigma g}).
\end{equation}
\end{linenomath}
We have two possible solutions, but it turns out that applying
the restriction $\mathbf{1}_D-p^\prime \ge 0$ selects the negative root.
To see this, we evaluate $\mathbf{1}_D-p^\prime$ at both roots:
\begin{linenomath}
\begin{equation}
\notag
\mathbf{1}_D - p^\prime = \frac{1}{2} (\mathbf{1}_D - p  \mp \sqrt{ (\mathbf{1}_D-p)^2 + 4 \sigma g}).
\end{equation}
\end{linenomath}
Using the fact that the data are non-negative, we have
\begin{linenomath}
\begin{equation}
\notag
\sqrt{ (\mathbf{1}_D-p)^2 + 4 \sigma g} \ge |\mathbf{1}_D-p|;
\end{equation}
\end{linenomath}
the positive root clearly leads to possible negative values for $\mathbf{1}_D-p^\prime$
while the negative root respects $\mathbf{1}_D-p^\prime \ge 0$ and yields Eq. (\ref{KLprox}).

For the final computation of a proximal mapping, we take a look at the data term
of the constrained, TV-minimization problem. From Eq. (\ref{ConstrainedTVFconj}), the
proximal mapping of interest is evaluated by:
\begin{linenomath}
\begin{equation}
\notag
prox_\sigma[F_1^*](p) = \argmin_{p^\prime}\; \left\{
\frac{\| p - p^\prime \|_2^2}{2 \sigma}
+\epsilon \| p^\prime \|_2 + \langle p^\prime,g \rangle_D \right\}.
\end{equation}
\end{linenomath}
Note the first term in the objective is spherically symmetric about $p$
and increasing with distance from $p$, and the
second term is also spherically symmetric about $\mathbf{0}_D$ and increasing
with distance from $\mathbf{0}_D$. If just these two terms
were present the minimum would lie on the line segment between $\mathbf{0}_D$ and $p$.
The third term, however, complicates the situation a little. We note that
this term is linear in $p^\prime$, and it can be combined with the first term by
completing the square. Performing this manipulation and ignoring constant
terms (independent of $p^\prime$) yields:
\begin{linenomath}
\begin{equation}
\notag
prox_\sigma[F_1^*](p) = \argmin_{p^\prime}\; \left\{
\frac{\| p - p^\prime - \sigma g \|_2^2}{2 \sigma}
+\epsilon \| p^\prime \|_2 \right\}.
\end{equation}
\end{linenomath}
By the geometric considerations discussed above, the minimizer lies
on the line segment between $\mathbf{0}_D$ and $p- \sigma g$. Analyzing this one-dimensional
minimization leads to Eq. (\ref{ConstrainedTVprox}).

\section{The finite differencing form of the image gradient and divergence}
\label{sec:graddiv}

In this appendix we write down the explicit forms of the finite differencing
approximations of $\nabla$ and $-div$ in two dimensions used in this article.
We use $x \in I$ to represent an $M \times M$ image and $x_{i,j}$ to refer to the $(i,j)$th pixel
of $x$. To specify the linear transform $\nabla$, we introduce the differencing
images $\Delta_s x \in I$ and  $\Delta_t x \in I$:
\begin{linenomath}
\begin{align}
\notag
\Delta_s x_{i,j} &=
\begin{cases}
x_{i+1,j} -x_{i,j} & i<M \\
-x_{i,j} & i=M
\end{cases}, \\
\notag
\Delta_t x_{i,j} &=
\begin{cases}
x_{i,j+1} -x_{i,j} & j<M \\
-x_{i,j} & j=M
\end{cases}.
\end{align}
\end{linenomath}
Using these definitions, $\nabla$ can be written as:
\begin{linenomath}
\begin{equation}
\notag
\nabla x = \left(
\begin{array}{c}
\Delta_s x \\
\Delta_t x
\end{array}
\right) .
\end{equation}
\end{linenomath}
With this form of $\nabla$, its transpose $-div$ becomes:
\begin{linenomath}
\begin{equation}
\notag
-div \left(
\begin{array}{c}
\Delta_s x \\
\Delta_t x
\end{array}
\right) = \left\{ -(\Delta_s x_{i,j} - \Delta_s x_{i-1,j}) -
 (\Delta_t x_{i,j} - \Delta_t x_{i,j-1}) , i \text{ and } j \in [1,M] \right\},
\end{equation}
\end{linenomath}
where the elements referred to outside the image border are set to zero:
$\Delta_s x_{0,j} =\Delta_s x_{i,0} =\Delta_t x_{0,j} =\Delta_t x_{i,0}=0$.
What the particular form of $\nabla$ is in its discrete form is not that important, but
it is critical that the discrete forms of $-div$ and $\nabla$ are the transposes of
each other.

\section{Preconditioned Chambolle-Pock algorithm demonstrated on the
KL-TV optimization problem}
\label{sec:preconditionedCP}

Chambolle and Pock followed their article, Ref. \cite{chambolle2011first},
with a pre-conditioned version of their algorithm that suits our
purpose of optimization problem proto-typing while potentially improving
algorithm efficiency substantially for the $\ell_2^2$-TV and KL-TV
optimization problems with small $\lambda$.  The new algorithm replaces the
constants $\sigma$ and $\tau$ with vector quantities that are computed
directly from the system matrix $K$, which yields a vector in space $Y$
from a vector in space $X$. One form of the suggested, diagonal pre-conditioners
uses the following weights:
\begin{linenomath}
\begin{align}
\label{precondsig}
\Sigma & =\frac{\mathbf{1}_Y}{|K| \mathbf{1}_X}, \\
\text{T} &= \frac{\mathbf{1}_X}{|K|^T \mathbf{1}_Y},
\end{align}
\end{linenomath}
where $\Sigma \in Y$, $\text{T} \in X$, and $|K|$ is the matrix formed
by taking the absolute value of each element of $K$. In order to generate
the CP algorithm instance incorporating pre-conditioning, the proximal mapping needs
to be modified:
\begin{linenomath}
\begin{equation}
\label{proximalPC}
prox_\Sigma[F](y) = \argmin_{y^\prime} \; \left\{ F(y^\prime) + \frac{1}{2}(y-y^\prime)^T
\left(\frac{ y - y^\prime}{\Sigma}\right) \right\}.
\end{equation}
\end{linenomath}
The second term in this minimization is still quadratic but no longer spherically
symmetric.
The difficulty in deriving the pre-conditioned CP algorithm instances is similar
to that of the original algorithm. On the one hand there is no need for finding
$\|K\|_2$, but on the other hand deriving the proximal mapping may become more
involved. For the $\ell_2^2$-TV and the KL-TV optimization problems,
the proximal mapping is simple to derive and it turns out that the mappings
can be arrived at by replacing $\sigma$ by $\Sigma$ and $\tau$ by $\text{T}$.

The gain in efficiency for small $\lambda$ comes from being able to absorb
this parameter into the TV term and allowing $\Sigma$ to account for the
mismatch between TV and data agreement terms. We modify the definitions of
$\nabla$ and $-div$ matrices from Appendix \ref{sec:graddiv}:
\begin{linenomath}
\begin{equation}
\notag
\nabla_\lambda x = \left(
\begin{array}{c}
\lambda \Delta_s x \\
\lambda \Delta_t x
\end{array}
\right) ,
\end{equation}
\end{linenomath}
and
\begin{linenomath}
\begin{equation}
\notag
-div_\lambda \left(
\begin{array}{c}
\Delta_s x \\
\Delta_t x
\end{array}
\right) = \left\{ -\lambda (\Delta_s x_{i,j} - \Delta_s x_{i-1,j}) -
\lambda (\Delta_t x_{i,j} - \Delta_t x_{i,j-1}) , i \text{ and } j \in [1,M] \right\},
\end{equation}
\end{linenomath}
where again the elements referred to outside the image border are set to zero:
$\Delta_s x_{0,j} =\Delta_s x_{i,0} =\Delta_t x_{0,j} =\Delta_t x_{i,0}=0$.

\begin{figure}[!h]
\begin{minipage}[b]{0.49\linewidth}
\centering
\centerline{\includegraphics[width=0.95\linewidth]{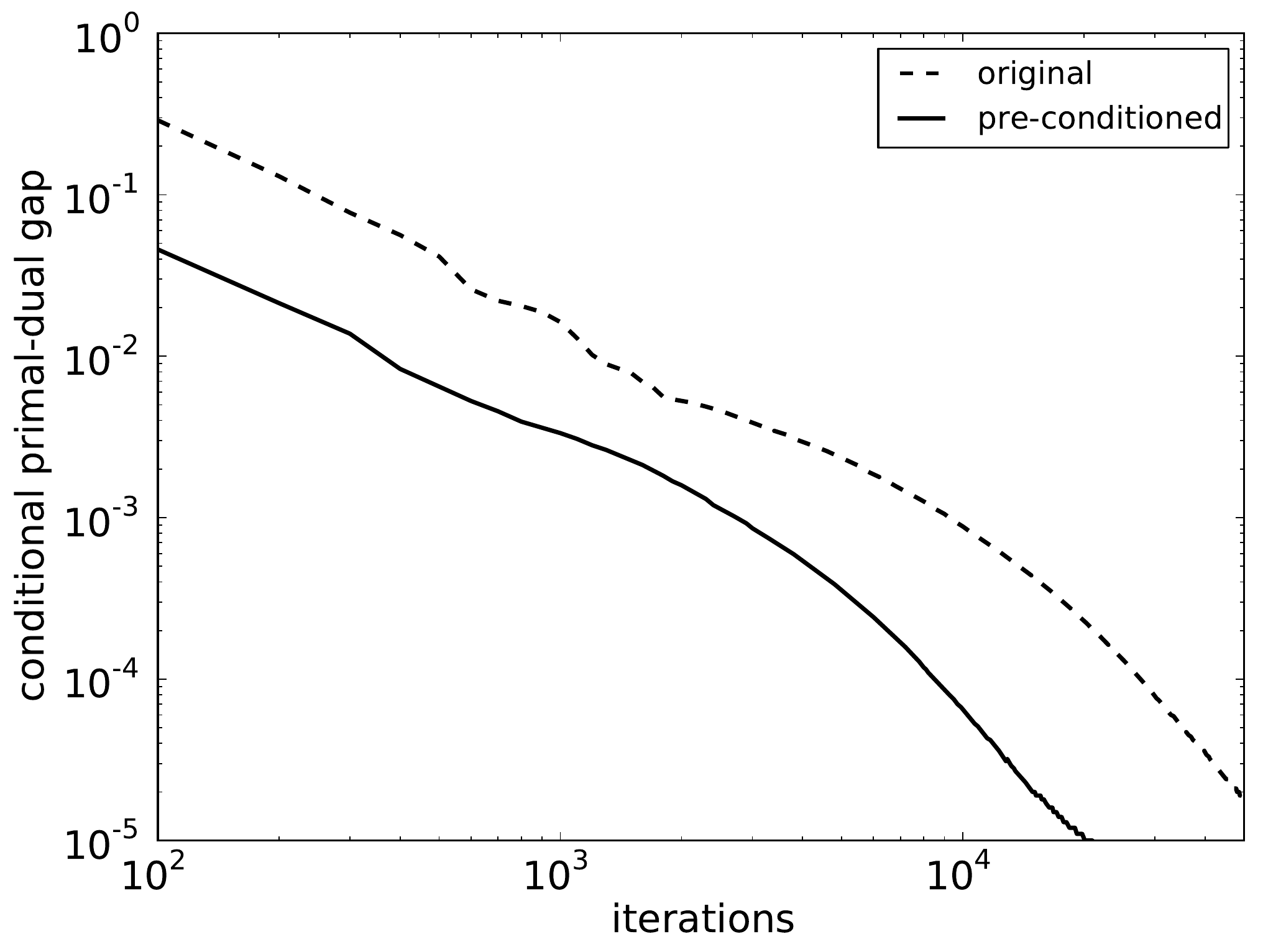}}
% \centerline{\includegraphics[width=8cm,clip=TRUE]{figs/viewscond.eps}}
\end{minipage}
\begin{minipage}[b]{0.49\linewidth}
\centering
\centerline{\includegraphics[width=0.95\linewidth]{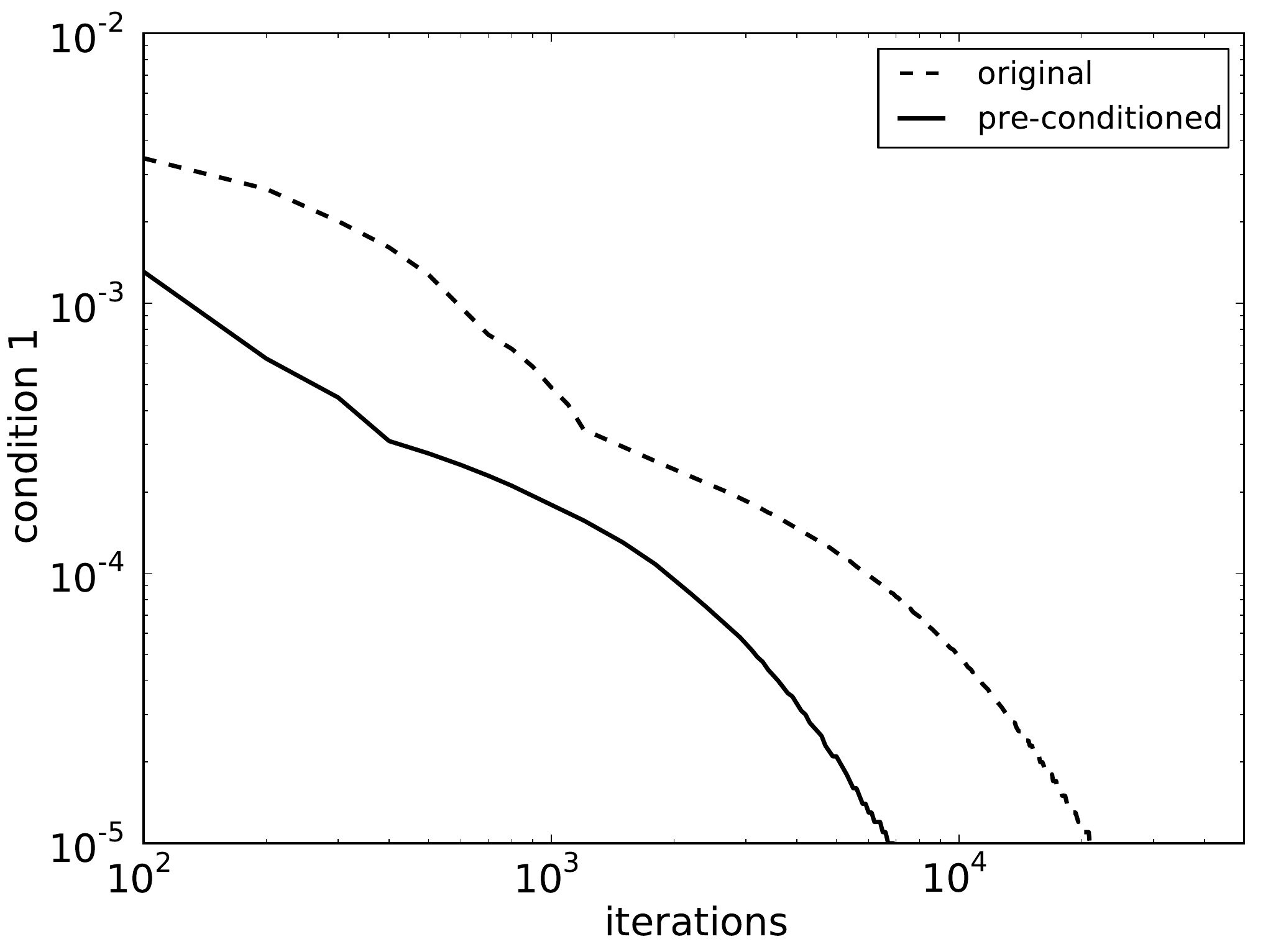}}
% \centerline{\includegraphics[width=8cm,clip=TRUE]{figs/viewscond.eps}}
\end{minipage}
\caption{
(Left) Convergence of the partial primal-dual gap for the CP algorithm instance solving 
Eq. (\ref{KLTVprimal}) for $\lambda=2 \times 10^{-5}$ for the original and
pre-conditioned CP algorithm.
(Right) Plot indicating agreement with condition 1 for the KL-TV optimization
problem. See Fig. \ref{fig:pdgap} for explanation.
\label{fig:PC}}
\end{figure}
For a complete example, we write the pre-conditioned CP algorithm instance for KL-TV 
in Listing \ref{algPCkltv}
\begin{algorithm}
\caption{Pseudocode for $N$-steps of the KL-TV pre-conditioned CP algorithm instance.}
\label{algPCkltv}
\begin{algorithmic}[1]
\STATE $\Sigma_1 \gets \mathbf{1}_D/(|A|\mathbf{1}_I); \;
\Sigma_2 \gets \mathbf{1}_V/(|\nabla_\lambda| \mathbf{1}_I); \;
\text{T}  \gets \mathbf{1}_I/(|A^T|\mathbf{1}_D \, + \, |div_\lambda| \mathbf{1}_V)$
\STATE $\theta \gets 1; \; n \gets 0$
\STATE initialize $u_0$, $p_0$, and $q_0$ to zero values
\STATE $\bar{u}_0 \gets u_0$
\REPEAT
\STATE $p_{n+1} \gets \frac{1}{2}
\left( \mathbf{1}_D + p_n+\Sigma_1 A\bar{u}_n - \sqrt{ ( p_n+\Sigma_1 A\bar{u}_n- \mathbf{1}_D)^2 + 4 \Sigma_1 g} \right)$
\label{pPCkltvupdate}
\STATE $q_{n+1} \gets (q_n + \Sigma_2 \nabla_\lambda \bar{u}_n )/
\max(\mathbf{1}_I,|q_n + \Sigma_2 \nabla_\lambda \bar{u}_n |)$ \label{qPCkltvupdate}
\STATE $u_{n+1} \gets  u_n - \text{T} A^T p_{n+1} + \text{T} div_\lambda \, q_{n+1}$ \label{uPCkltvupdate}
\STATE $\bar{u}_{n+1} \gets u_{n+1} + \theta(u_{n+1} - u_n)$
\STATE $n \gets n+1$
\UNTIL{$n \ge N$}
\end{algorithmic}
\end{algorithm}
To illustrate the potential gain in efficiency, we show the condition
primal-dual gap as a function of iteration number for the KL-TV problem
with $\lambda=2 \times 10^{-5}$ in Fig. \ref{fig:PC}. 
While we have presented the pre-conditioned CP algorithm as a patch
for the small $\lambda$ case, it really provides an alternative prototyping
algorithm and it can be used instead of the original CP algorithm. \\
%Phew! that was a lot of equations.

\newpage

\bibliographystyle{IEEEtran}
\bibliography{sampling}

\begin{thebibliography}{10}
\providecommand{\url}[1]{#1}
\csname url@samestyle\endcsname
\providecommand{\newblock}{\relax}
\providecommand{\bibinfo}[2]{#2}
\providecommand{\BIBentrySTDinterwordspacing}{\spaceskip=0pt\relax}
\providecommand{\BIBentryALTinterwordstretchfactor}{4}
\providecommand{\BIBentryALTinterwordspacing}{\spaceskip=\fontdimen2\font plus
\BIBentryALTinterwordstretchfactor\fontdimen3\font minus
  \fontdimen4\font\relax}
\providecommand{\BIBforeignlanguage}[2]{{%
\expandafter\ifx\csname l@#1\endcsname\relax
\typeout{** WARNING: IEEEtran.bst: No hyphenation pattern has been}%
\typeout{** loaded for the language `#1'. Using the pattern for}%
\typeout{** the default language instead.}%
\else
\language=\csname l@#1\endcsname
\fi
#2}}
\providecommand{\BIBdecl}{\relax}
\BIBdecl

\bibitem{McCollough:09}
C.~H. McCollough, A.~N. Primak, N.~Braun, J.~Kofler, L.~Yu, and J.~Christner,
  ``Strategies for reducing radiation dose in {CT},'' \emph{Radiol. Clin. N.
  Am.}, vol.~47, pp. 27--40, 2009.

\bibitem{PanIP:09}
X.~Pan, E.~Y. Sidky, and M.~Vannier, ``Why do commercial {CT} scanners still
  employ traditional, filtered back-projection for image reconstruction?''
  \emph{Inv. Prob.}, vol.~25, pp. 123\,009--(1--36), 2009.

\bibitem{ziegler2008iterative}
A.~Ziegler, T.~Nielsen, and M.~Grass, ``Iterative reconstruction of a region of
  interest for transmission tomography,'' \emph{Med. Phys.}, vol.~35, pp.
  1317--1327, 2008.

\bibitem{li2002accurate}
M.~Li, H.~Yang, and H.~Kudo, ``{An accurate iterative reconstruction algorithm
  for sparse objects: Application to 3D blood vessel reconstruction from a
  limited number of projections},'' \emph{Phy. Med. Biol.}, vol.~47, pp.
  2599--2609, 2002.

\bibitem{SidkyTV:06}
E.~Y. Sidky, C.-M. Kao, and X.~Pan, ``Accurate image reconstruction from
  few-views and limited-angle data in divergent-beam {CT},'' \emph{J. X-ray
  Sci. Tech.}, vol.~14, pp. 119--139, 2006.

\bibitem{Sidky2008image}
E.~Y. Sidky and X.~Pan, ``{Image reconstruction in circular cone-beam computed
  tomography by constrained, total-variation minimization},'' \emph{Phys. Med.
  Biol.}, vol.~53, pp. 4777--4807, 2008.

\bibitem{Chen2008prior}
G.~H. Chen, J.~Tang, and S.~Leng, ``Prior image constrained compressed sensing
  ({PICCS}): a method to accurately reconstruct dynamic {CT} images from highly
  undersampled projection data sets,'' \emph{Med. Phys.}, vol.~35, pp.
  660--663, 2008.

\bibitem{SidkyPC:10}
E.~Y. Sidky, M.~A. Anastasio, and X.~Pan, ``Image reconstruction exploiting
  object sparsity in boundary-enhanced x-ray phase-contrast tomography,''
  \emph{Opt. Express}, vol.~18, pp. 10\,404--10\,422, 2010.

\bibitem{ritschl2011improved}
L.~Ritschl, F.~Bergner, C.~Fleischmann, and M.~Kachelrie{\ss}, ``Improved total
  variation-based {CT} image reconstruction applied to clinical data,''
  \emph{Phys. Med. Biol.}, vol.~56, pp. 1545--1562, 2011.

\bibitem{Defrise:11}
M.~Defrise, C.~Vanhove, and X.~Liu, ``An algorithm for total variation
  regularization in high-dimensional linear problems,'' \emph{Inv. Prob.},
  vol.~27, p. 065002, 2011.

\bibitem{Fessler:2011}
S.~Ramani and J.~Fessler, ``A splitting-based iterative algorithm for
  accelerated statistical {X}-ray {CT} reconstruction,'' \emph{IEEE Trans. Med.
  Imag.}, 2011, available online at IEEE TMI - early access.

\bibitem{jakob:2011}
J.~H. J{\o}rgensen, E.~Y. Sidky, and X.~Pan, ``Analysis of discrete-to-discrete
  imaging models for iterative tomographic image reconstruction and compressive
  sensing,'' 2011, arxiv preprint arxiv:1109.0629
  (http://arxiv.org/abs/1109.0629).

\bibitem{Nocedal:06}
J.~Nocedal and S.~Wright, \emph{Numerical Optimization, 2nd ed.}\hskip 1em plus
  0.5em minus 0.4em\relax Springer, 2006.

\bibitem{erdogan1999ordered}
H.~Erdogan and J.~A. Fessler, ``Ordered subsets algorithms for transmission
  tomography,'' \emph{Phys. Med. Biol}, vol.~44, pp. 2835--2852, 1999.

\bibitem{green1984iteratively}
P.~Green, ``Iteratively reweighted least squares for maximum likelihood
  estimation, and some robust and resistant alternatives,'' \emph{J. Royal
  Stat. Soc., Ser. B}, vol.~46, pp. 149--192, 1984.

\bibitem{yin2008bregman}
W.~Yin, S.~Osher, D.~Goldfarb, and J.~Darbon, ``Bregman iterative algorithms
  for $\ell$1-minimization with applications to compressed sensing,''
  \emph{SIAM J. Imag. Sci.}, vol.~1, pp. 143--168, 2008.

\bibitem{combettes2008proximal}
P.~L. Combettes and J.~C. Pesquet, ``A proximal decomposition method for
  solving convex variational inverse problems,'' \emph{Inv. Prob.}, vol.~24,
  pp. 065\,014--27pp, 2008.

\bibitem{Beck:09}
A.~Beck and M.~Teboulle, ``Fast gradient-based algorithms for constrained total
  variation image denoising and deblurring problems,'' \emph{IEEE Trans. Imag.
  Proc.}, vol.~18, pp. 2419--2434, 2009.

\bibitem{becker2010templates}
S.~R. Becker, E.~J. Candes, and M.~Grant, ``Templates for convex cone problems
  with applications to sparse signal recovery,'' 2010, arxiv preprint
  arXiv:1009.2065.

\bibitem{chambolle2011first}
A.~Chambolle and T.~Pock, ``A first-order primal-dual algorithm for convex
  problems with applications to imaging,'' \emph{J. Math. Imag. Vis.}, vol.~40,
  pp. 1--26, 2011.

\bibitem{jakob:11}
T.~L. Jensen, J.~H. J{\o}rgensen, P.~C. Hansen, and S.~H. Jensen,
  ``Implementation of an optimal first-order method for strongly convex total
  variation regularization,'' 2011, accepted to {BIT}, available at
  http://arxiv.org/abs/1105.3723.

\bibitem{Choi:10}
K.~Choi, J.~Wang, L.~Zhu, T.-S. Suh, S.~Boyd, and L.~Xing, ``Compressed sensing
  based cone-beam computed tomography reconstruction with a first-order
  method,'' \emph{Med. Phys.}, vol.~37, pp. 5113--5125, 2010.

\bibitem{joergensen2011toward}
J.~H. J{\o}rgensen, P.~C. Hansen, E.~Y. Sidky, I.~S. Reiser, and X.~Pan,
  ``Toward optimal {X}-ray flux utilization in breast {CT},'' in
  \emph{Proceedings of the 11th International Meeting on Three-Dimensional
  Image Reconstruction in Radiology and Nuclear Medicine}, 2011, arxiv preprint
  arxiv:1104.1588 (http://arxiv.org/abs/1104.1588).

\bibitem{rockafellar1970convex}
R.~T. Rockafellar, \emph{Convex analysis}.\hskip 1em plus 0.5em minus
  0.4em\relax Princeton Univ. Press, 1970.

\bibitem{Barrett:FIS}
H.~H. Barrett and K.~J. Myers, \emph{Foundations of Image Science}.\hskip 1em
  plus 0.5em minus 0.4em\relax Hoboken, NJ: John Wiley \& Sons, 2004.

\bibitem{Siddon:1985}
R.~L. Siddon, ``Fast calculation of the exact radiological path for a
  three-dimensional {CT} array,'' \emph{Med. Phys.}, vol.~12, pp. 252--255,
  1985.

\bibitem{de2004distance}
B.~D. Man and S.~Basu, ``Distance-driven projection and backprojection in three
  dimensions,'' \emph{Phys. Med. Biol}, vol.~49, pp. 2463--2475, 2004.

\bibitem{Mueller:07}
F.~Xu and K.~Mueller, ``Real-time 3{D} computed tomographic reconstruction
  using commodity graphics hardware,'' \emph{Phys. Med. Biol.}, vol.~52, pp.
  3405--3419, 2007.

\bibitem{zeng2000unmatched}
G.~L. Zeng and G.~T. Gullberg, ``Unmatched projector/backprojector pairs in an
  iterative reconstruction algorithm,'' \emph{IEEE Trans. Med. Imag.}, vol.~19,
  pp. 548--555, 2000.

\bibitem{Bian:10}
J.~Bian, J.~H. Siewerdsen, X.~Han, E.~Y. Sidky, J.~L. Prince, C.~A. Pelizzari,
  and X.~Pan, ``Evaluation of sparse-view reconstruction from
  flat-panel-detector cone-beam {CT},'' \emph{Phys. Med. Biol}, vol.~55, pp.
  6575--6599, 2010.

\bibitem{han2011algorithm}
X.~Han, J.~Bian, D.~R. Eaker, T.~L. Kline, E.~Y. Sidky, E.~L. Ritman, and
  X.~Pan, ``Algorithm-enabled low-dose micro-{CT} imaging,'' \emph{IEEE Trans.
  Med. Imag.}, vol.~30, pp. 606--620, 2011.

\bibitem{xia:043706}
D.~Xia, X.~Xiao, J.~Bian, X.~Han, E.~Y. Sidky, F.~D. Carlo, and X.~Pan, ``Image
  reconstruction from sparse data in synchrotron-radiation-based
  microtomography,'' \emph{Rev. Sci. Inst.}, vol.~82, pp. 043\,706 -- 9 pgs.,
  2011.

\bibitem{sidky2011special}
E.~Y. Sidky, Y.~Duchin, C.~Ullberg, and X.~Pan, ``{X}-ray computed tomography:
  advances in image formation: {A} constrained, total-variation minimization
  algorithm for low-intensity {X}-ray {CT},'' \emph{Med. Phys.}, vol.~38, pp.
  S117--S125, 2011.

\bibitem{elad2010sparse}
M.~Elad, \emph{Sparse and redundant representations: from theory to
  applications in signal and image processing}.\hskip 1em plus 0.5em minus
  0.4em\relax Springer Verlag, 2010.

\bibitem{chan2005aspects}
T.~F. Chan and S.~Esedo{$\bar{\text{g}}$}lu, ``Aspects of total variation
  regularized {$L^1$} function approximation,'' \emph{SIAM J. Appl. Math.},
  vol.~65, pp. 1817--1837, 2005.

\bibitem{boyd2004convex}
S.~P. Boyd and L.~Vandenberghe, \emph{Convex optimization}.\hskip 1em plus
  0.5em minus 0.4em\relax Cambridge University Press, 2004.

\bibitem{candes2008introduction}
E.~J. Cand{\`e}s and M.~B. Wakin, ``An introduction to compressive sampling,''
  \emph{IEEE Sig. Proc. Mag.}, vol.~25, pp. 21--30, 2008.

\bibitem{reiser2010task}
I.~Reiser and R.~M. Nishikawa, ``Task-based assessment of breast tomosynthesis:
  Effect of acquisition parameters and quantum noise,'' \emph{Med. Phys.},
  vol.~37, pp. 1591--1600, 2010.

\bibitem{Pock2011}
T.~Pock and A.~Chambolle, ``Diagonal preconditioning for first order
  primal-dual algorithms in convex optimization,'' in \emph{International
  Conference on Computer Vision ({ICCV} 2011)}, 2011, to appear; available at
  URL: http://gpu4vision.icg.tugraz.at.

\end{thebibliography}

\end{document}